\documentclass[12pt]{amsart}
\usepackage{amsfonts,amssymb,float,textcomp}
\RequirePackage{ifpdf}
\usepackage{amsmath}
\usepackage{color}
\usepackage{pstricks}
\usepackage{pst-node}
\usepackage{pst-plot}
\usepackage{pst-text}
\usepackage{pst-slpe}
\usepackage{pstricks-add}
\usepackage{tikz-cd}

\def\appendix#1{
\addtocounter{section}{1} \setcounter{equation}{0}
\renewcommand{\thesection}{\Alph{section}}
\section*{Appendix \thesection\protect\indent\quad
#1}
}
\renewcommand{\theequation}{\thesection.\arabic{equation}}

\def\hybrid{\topmargin 0pt      \oddsidemargin 0pt
        \headheight 0pt \headsep 0pt
       \textwidth 6.5in        
        \textheight 9.0in       
        \marginparwidth 0.0in
        \parskip 5pt plus 1pt   \jot = 1.5ex}
\catcode`\@=11
\def\marginnote#1{}

\newcount\hour
\newcount\minute
\newtoks\amorpm
\hour=\time\divide\hour by60 \minute=\time{\multiply\hour by60
\global\advance\minute by-\hour}
\edef\standardtime{{\ifnum\hour<12 \global\amorpm={am}%
        \else\global\amorpm={pm}\advance\hour by-12 \fi
        \ifnum\hour=0 \hour=12 \fi
        \number\hour:\ifnum\minute<10 0\fi\number\minute\the\amorpm}}
\edef\militarytime{\number\hour:\ifnum\minute<100\fi\number\minute}

\newcommand{\tcy}{\textcolor{yellow}}
\newcommand{\tcr}{\textcolor{red}}
\newcommand{\tcb}{\textcolor{blue}}
\newcommand{\tcg}{\textcolor{green}}
\newcommand{\tcw}{\textcolor{white}}
\newcommand{\mapi}{{\stackrel{{i}}{{\to}}}}
\let\hhat=\widehat
\let\ttilde=\widetilde

\def\draftlabel#1{{\@bsphack\if@filesw {\let\thepage\relax
      \xdef\@gtempa{\write\@auxout{\string
          \newlabel{#1}{{\@currentlabel}{\thepage}}}}}\@gtempa \if@nobreak
    \ifvmode\nobreak\fi\fi\fi\@esphack} \gdef\@eqnlabel{#1}}
    \def\@eqnlabel{}
\def\@vacuum{}
\def\draftmarginnote#1{\marginpar{\raggedright\scriptsize\tt#1}}

\def\draft{
  \oddsidemargin -.5truein
  \def\@oddfoot{\footnotesize \sl preliminary draft \hfil
    \rm\thepage\hfil\sl\today\quad\militarytime}
  \let\@evenfoot\@oddfoot \overfullrule 3pt
    \let\label=\draftlabel
    \let\marginnote=\draftmarginnote
  \def\@eqnnum{(\theequation)\rlap{\kern\marginparsep\tt\@eqnlabel}%
    \global\let\@eqnlabel\@vacuum}

  }

\voffset= - 0.5in \hoffset= - 1.0in

\newcommand{\tr}{\,{\rm Tr}\,}
\newcommand{\ID}{1\!\!1}
\def\d{\partial}
\def\eps{\varepsilon}
\def\be{\begin{equation}}
\def\ee{\end{equation}}
\def\bea{\begin{eqnarray}}
\def\eea{\end{eqnarray}}
\def\<{\langle}
\def\>{\rangle}
\def\nn{\nonumber}
\def\HH{{\Bbb H}}
\def\RR{{\Bbb R}}
\def\QQ{{\Bbb Q}}
\def\CC{{\Bbb C}}
\def\ZZ{{\Bbb Z}}
\def\A{\mathcal{A}}
\def\B{\mathcal{B}}
\def\X{\mathcal{X}}
\def\FF{\mathfrak{F}}
\def\bb{{\bf b}}
\def\M{{\mathcal M}}
\let\wtd=\widetilde
\def\Tr{{\rm Tr}}
\def\g{{\mathfrak g}}
\def\gl{{\mathfrak{gl}}}
\def\slf{{\mathfrak{sl}}}
\def\G{{\mathcal G}}
\def\h{\mathfrak h}
\def\D{{\mathfrak d}}
\def\U{{\mathcal U}}
\def\w{{\mathbf w}}
\def\u{{\mathbf u}}
\def\v{{\mathbf v}}
\def\vv{{\mathbf v}}
\def\otim{\mathop{\otimes}}
\def\tr{\mathop{\rm{tr}}}
\def\thi{\vartheta}
\def\ocomma{{\phantom{\Bigm|}^{\phantom {X}}_{\raise-1.5pt\hbox{,}}\!\!\!\!\!\!\otimes}}
\def\Mat{\operatorname{Mat}}
\def\End{\operatorname{End}}
\def\Id{{\operatorname {Id}}}

\def\Poi{{\{\cdot,\cdot\}}}
\newcommand{\col}[1]{{\raise-2pt\hbox{\tiny$\bullet$}\hskip -4.5pt \raise4pt\hbox{\tiny$\bullet$}{{#1}} \raise-2pt\hbox{\tiny$\bullet$}\hskip -4.5pt \raise4pt\hbox{\tiny$\bullet$}}}

\newcommand{\sheet}[2]{{\stackrel{{#1}}{{#2}}}}
\newcommand{\groupoid}{{\stackrel{{\rightarrow}}{{_\rightarrow}}}}
\def\dnabla{{\raisebox{2pt}{$\bigtriangledown$}}\negthinspace}

\newtheorem{theorem}{Theorem}[section]
\newtheorem{lemma}[theorem]{Lemma}

\theoremstyle{definition}
\newtheorem{definition}[theorem]{Definition}

\newtheorem{remark}[theorem]{Remark}

\newtheorem{corollary}[theorem]{Corollary}

\allowdisplaybreaks

\long\def\rem#1{}

\begin{document}

\title[Characteristic equation for symplectic groupoid and cluster algebras]
{Characteristic equation for symplectic groupoid and cluster algebras}
\author{Leonid O. Chekhov$^{\ast}$}
\thanks{$^{\ast}$ Steklov Mathematical
Institute, Moscow, Russia, National Research University Higher School of Economics, Russia, and Michigan State University, East Lansing, USA. Email: chekhov@msu.edu.} 
\author{Michael Shapiro$^{\ast\ast}$} 
\thanks{$^{\ast\ast}$ Michigan State University, East Lansing, USA and National Research University Higher School of Economics, Russia. Email: mshapiro@msu.edu.} 
\author{Huang Shibo$^{\ast\ast\ast}$}
\thanks{$^{\ast\ast\ast}$  Xi'an Jiaotong University, Xi'an, Shaanxi, P.R. China.}

\begin{abstract}
We use the Darboux coordinate representation found by two of the authors (L.Ch. and M.Sh.) for entries of general symplectic leaves of the $\mathcal A_n$-groupoid of upper-triangular matrices to express roots of the characteristic equation $\det(\mathbb A-\lambda \mathbb A^{\text{T}})=0$, with $\mathbb A\in \mathcal A_n$, in terms of Casimirs of this Darboux coordinate representation, which is based on cluster variables of Fock--Goncharov higher Teichm\"uller spaces for the algebra $sl_n$. We show that roots of the characteristic equation are simple monomials of cluster Casimir elements. This statement remains valid in the quantum case as well. We consider a generalization of $\mathcal A_n$-groupoid to a $\mathcal A_{Sp_{2m}}$-groupoid.
\end{abstract}

\maketitle

\section{Introduction}\label{s:intro}
\setcounter{equation}{0}

\subsection{Symplectic groupoid and induced Poisson structure on the unipotent upper triangular matrices}

Let $\mathcal A_n\subseteq gl_n$ be a subspace of unipotent upper-triangular $n\times n$ matrices. 
We identify elements $\mathbb A$ of $\mathcal A_n$ with matrices of bilinear forms on $\mathbb{C}^n$. 
 The matrix $B\in GL_n$ of a change of a basis  in $\mathbb{C}^n$ 
 takes a matrix of bilinear form  $\mathbb A\in \mathcal A_n$ to $B\mathbb AB^{\text{T}}$ determining a {\sl dynamics of form transformations}. A most interesting case is when the transformed form $B\mathbb AB^{\text{T}}$ lies in the same class $\mathcal A_n$. For a fixed  $\mathbb A\in\mathcal A_n$ all $B\in GL_n$ such that $B{\mathbb A} B^T\in{\mathcal A}_n$ form  a submanifold of  $GL_n$, which is not a subspace and may have a very involved structure.  We introduce the space of \emph{morphisms} identified with admissible pairs of matrices $(B,\mathbb A)$ such that  
\be
\label{eq:tr}
\mathcal M=\bigl\{ (B,\mathbb A) \bigm | B\in GL_n 
,\ \mathbb A\in \mathcal A_n,\ B\mathbb AB^{\text{T}}\in \mathcal A_n  \bigr\}.
\ee

In 2000, Bondal \cite{Bondal} obtained the Poisson structure on $\mathcal A_n$ using the algebroid construction: For $B=e^{\mathfrak g}$, we first define the anchor map $D_{\mathbb A}$ to the tangent space $ T \mathcal A_n$:
\be
\label{gA}
\begin{array}{lccc}
D_{\mathbb A}:&\mathfrak g_{\mathbb A}&\to& T \mathcal A_n\\
& g &\mapsto &{\mathbb A}g+g^{\text{T}}{\mathbb A}.\\ \end{array}
\ee
where  $\mathfrak{g}_{\mathbb A}$ is the linear subspace
$$
\mathfrak{g}_{\mathbb A}:=\left\{ g\in\gl_{n}(\mathbb C),|\,
\mathbb A+ {\mathbb A}g+g^{\text{T}}{\mathbb A}\in\mathcal A_n \right\}
$$
of elements $g$ leaving $\mathbb A$ unipotent.

\begin{lemma}\label{lemma-g}\cite{Bondal}
The map
\be
\label{P_A}
\begin{array}{lccl}
P_{\mathbb A}:&T^\ast \mathcal A_n&\to& \mathfrak g_{\mathbb A}\\
& w &\mapsto &P_{-,1/2}(w{\mathbb A})-P_{+,1/2}(w^{\text{T}}{\mathbb A}^{\text{T}}),\\ \end{array}
\ee
where $P_{\pm,1/2}$ are the projection operators:
\be\label{eq:pr}
P_{\pm,1/2}a_{i,j}:=\frac{1\pm {\rm sign}(j-i)}{2}a_{i,j}, \quad i, j=1,\dots,n,
\ee
and $w\in T^\ast\mathcal A_n$ is a strictly lower triangular matrix, 
defines an isomorphism between the Lie algebroid $(\mathfrak g,D_{\mathbb A})$ and the Lie algebroid
$\left(T^\ast \mathcal A_n,D_{\mathbb A}P_{\mathbb A}  \right)$.
\end{lemma}

The Poisson bi-vector $\Pi$ on $\mathcal A_n$ is then obtained by the anchor map on the Lie algebroid
$\left(T^\ast \mathcal A_n,D_{\mathbb A}P_{\mathbb A}  \right)$ (see Proposition 10.1.4 in \cite{kirill}) as:
\be
\label{eq:biv}
\begin{array}{lccl}
\Pi:&T^\ast \mathcal A_n\times T^\ast \mathcal A_n
&\mapsto& \mathcal C^\infty(\mathcal A)\\
&(\omega_1,\omega_2)&\to&\Tr\left(\omega_1 D_{\mathbb A}P_{\mathbb A}  (\omega_2)
\right)
\end{array}
\ee
It can be checked explicitly that the above bilinear form is in fact skew-symmetric and gives rise to the Poisson bracket
\be
\label{Poisson-bracket}
\{a_{i,k},a_{j,l}\}:=\frac{\partial}{\partial{\rm d}a_{i,k}}\wedge\frac{\partial}{\partial{\rm d}a_{j,l}}
\Tr\left( {\rm d}a_{i,k} D_{\mathbb A}P_{\mathbb A}  ({\rm d}a_{j,l})\right),
\ee
having the following form: 
\begin{eqnarray}\label{eq:du}
&&
\left\{a_{ik},a_{jl}\right\}=0,\quad\hbox{for}\ i<k<j<l,\hbox{ and }  i<j<l<k,\nn
\\&&
\left\{a_{ik},a_{jl}\right\}=2 \left(a_{ij}a_{kl}-a_{il}a_{kj}\right),\quad\hbox{for}\  i<j<k<l,
\\&&
\left\{a_{ik},a_{kl}\right\}=a_{ik}a_{kl}-2a_{il},\quad\hbox{for}\  i<k<l,\nn
\\&&
\left\{a_{ik},a_{jk}\right\}=-a_{ik}a_{jk}+2a_{ij},\quad\hbox{for} \  i< j<k,\nn
\\&&
\left\{a_{ik},a_{il}\right\}=-a_{ik}a_{il}+a_{kl},\quad\hbox{for} \  i<k<l.\nn
\end{eqnarray}
This bracket is famous in mathematical physics and is known as Gavrilik--Klimyk--Nelson--Regge--Dubrovin--Ugaglia bracket \cite{GK91,NR,NRZ,Dub,Ugaglia}. It originally arose in the context of 2D quantum gravity; technically it governs from skein relations satisfied by a specially chosen \cite{ChF3} finite subset of {\it geodesic functions} $a_{ik}$ (traces of monodromies of $SL_2$ Fuchsian systems), which are principal observables in the theory of 2D gravity; a log-canonical (Darboux) bracket on the space of Thurston shear coordinates $z_\alpha$  on the Teichm\"uller space $T_{g,s}$ of Riemann surfaces $\Sigma_{g,s}$ of genus $g$ with $s=1,2$ holes was shown \cite{ChF2} to induce the above bracket on $a_{ik}$ for $n=2g+s$. All geodesic functions on any  $\Sigma_{g,s}$ admit an explicit combinatorial description \cite{F97} and they are Laurent polynomials with positive integer coefficients of $e^{z_\alpha/2}$. The Poisson bracket of $z_\alpha$ spanning the Teichm\"uller space $T_{g,s}$ has exactly $s$ Casimirs, which are linear combinations of shear coordinates incident to the holes, so the subspace of $z_{\alpha}$ orthogonal to the subspace of Casimirs parameterizes a symplectic leaf in the Teichm\"uller space called a {\it geometric symplectic leaf}.

However the Poisson dimensions of $T_{g,s}$ $(6g-6+s)$ and those of $\mathcal A_n$ $(n(n-1)/2-[n/2])$ become different for $n\ge 6$, so a novel Darboux coordinate construction for higher-dimensional symplectic leaves of $\mathcal A_n$ was found in \cite{ChSh2} where elements of $\mathcal A_n$ were parameterized by cluster variables of the Fock--Goncharov framed moduli space $\X_{SL_n,\Sigma}$, where $\Sigma$ is a disk with three marked points
\cite{FG1}.

Note that algebra (\ref{eq:du}) is universal: any choice of the subspace $\mathcal A\subseteq sl_n$ for a $(B,\mathbb A)$-system subject to the form dynamics results in Poisson relations  (\ref{eq:du}) provided $\mathcal A$ is a Lagrangian submanifold \cite{ChM3}, so constructing cluster realizations of $\mathbb A$ outside the unipotent upper-triangular case is also important. We present this construction for the first natural extension of $\mathcal A_n$, which is $\mathcal A_{Sp_{2m}}$.  The quantum analogue of (\ref{eq:du}) is known as the {\sl quantum reflection equation} (see Theorem~\ref{th:A}), and the cluster coordinate realization of a general symplectic leaf of the quantum reflection equation was also constructed in \cite{ChSh2}.

An important problem is to construct Casimirs of algebras (\ref{eq:du}) in semiclassical and quantum cases. In the semiclassical limit, Bondal showed that the algebra $\mathcal A_n$ admits $[n/2]$ algebraically independent Casimir elements, which are coefficients of the reciprocal polynomial $\det [\mathbb A-\lambda\mathbb A^{\text{T}}]$; in the geometric cases $n=3$ and $n=4$, these coefficients turn out to be simple functions of {\em hole perimeters} $p_i=\sum z_{\alpha_i}$, which are always linear functions of shear coordinates. This statement remains valid in the quantum case as well. At the same time, being expressed in $a_{ij}$, the same coefficients become inhomogeneous functions  producing the quantum Markov invariant for $n=3$ and quantum invariants first obtained by Bullock and Przytycki ``by brute force'' in \cite{BP99} for $n=4$.

Although we do not exploit it in this text, note that the space $\mathcal A_n$ admits a discrete {\em braid-group action} generated by special morphisms $\beta_{i,i+1}:{\mathcal A_n}\to{\mathcal A_n}$, $i=1,\dots,n-1$, such that $\beta_{i,i+1}[{\mathbb A}]:=B_{i,i+1}{\mathbb A}B_{i,i+1}^T\in \mathcal A_n$; these transformations correspond to Dehn twists in the geometrical case, and matrices $B_{i,i+1}$ differ from the unity matrix only in their diagonal $2\times 2$-blocks having the structure
$$
\begin{array}{c|cc}
&i &i+1\\
\hline
i& q^{1/2}a^\hbar_{i,i+1}&-q\phantom{\Bigm|}\\
i+1&1&0
\end{array}
$$
in the quantum case constructed in  \cite{ChM}; the matrix $\mathbb A$ then becomes a 
{\sl quantum} matrix $\mathbb A^\hbar$:
\be\label{A-quantum}
\mathbb A^\hbar:=\left[ \begin{array}{ccccc}
q^{-1/2}  & a^\hbar_{1,2}  & a^\hbar_{1,3}  &\dots & a^\hbar_{1,n}   \\
 0 & q^{-1/2}  & a^\hbar_{2,3} & \dots  & a^\hbar_{2,n}  \\
 0 & 0 & q^{-1/2}  & \ddots &   \vdots \\
 \vdots & \vdots   & \ddots & \ddots &  a^\hbar_{n-1,n} \\
 0 & 0 & \dots & 0 & q^{-1/2} 
 \end{array}\right],
\ee
where $a^\hbar_{i,j}$ are self-adjoint $\bigl( \bigl[a^\hbar_{i,j}\bigr]^\star=a^\hbar_{i,j}\bigr)$ operators enjoying quadratic--linear algebraic relations following from the  quantum reflection equation (see Theorem~\ref{th:A}). The quantum braid-group action is $\beta^{\hbar}_{i,i+1} : \mathbb A^\hbar\to B^\hbar_{i,i+1}\mathbb A^\hbar \bigl[ B^\hbar_{i,i+1} \bigr]^\dagger$ with $\dagger$ standing for the Hermitian conjugation. In  \cite{ChSh2}, the braid-group action was realized as sequences of cluster mutations in the $\mathcal A_n$-quiver  obtained by a special amalgamation procedure (described below) from the quiver for $\X_{SL_n,\Sigma}$.


Our goal in this paper is to study the structure of resolvents of the quantum operators $\mathbb A\bigl[\mathbb A^{\dagger}\bigr]^{-1}:=\mathbb A\mathbb A^{-\dagger}$. This combination plays a prominent role: first, assuming all operatorial expressions be invertible, for $\mathbb A$ transforming as a form (\ref{eq:tr}), this combination undergoes an adjoint transformation,
\be
\mathbb A\mathbb A^{-\dagger}\to B \mathbb A\mathbb A^{-\dagger} B^{-1},
\ee
so we can address the problem of solving the resolvent equation: for which $\lambda\in\mathbb C$ the operator $\mathbb A\mathbb A^{-\dagger}-\lambda\mathbf I$ admits null vectors? We completely solve this problem in Theorem~\ref{th:M} expressing all $\lambda$ in question via Casimirs of the $\mathcal A_n$-quiver. 

With this solution in hands, we can address several intriguing problems: note that in the semiclassical case, the Jordan form of  $\mathbb A\mathbb A^{-\dagger}$ determines the dimension of the corresponding Poisson leaf. Although our technique does not allow keeping control over the Jordan form structure in case of coinciding eigenvalues $\lambda_i$, we can fully treat the case where all $n$ eigenvalues $\lambda_i$ are different. Then the Jordan form is diagonal and, assuming the existence of natural $N\ge n$ for which $\lambda_i^N=1$ for all $i$, we obtain that the $N$th power of the Jordan form is a unit matrix, and therefore $\bigl[\mathbb A\mathbb A^{-\dagger}\bigr]^N=\mathbf I$. This relation opens a way to constructing minimal-model  representations (containing null vectors in the sense of Kac) of related non-Abelian Verma-type modules generated by the action of $a_{i,j}$ on a proper vacuum vector, where the action of canonical cluster variables plays a role of a Dotsenko--Fatteev free-field representation. We may also address finite-dimensional reductions of affine Lie--Poisson systems.

\section{$b_n$-Algebras for the triangle $\Sigma_{0,1,3}$}\label{s:triangle}
\setcounter{equation}{0}

Let $\Sigma_{g,s,p}$ denote a topological genus $g$ surface with $s$ boundary components and $p>0$ marked points on the boundaries. Our main example is the disk with 3 marked points on the boundary, $\Sigma_{0,1,3}$, which corresponds to an ideal triangle in an ideal-triangle decomposition of a general $\Sigma_{g,s,p}$ (which comprises exactly $4g-4+2s+p$ ideal triangles). 

\subsection{Notations} 
Let lattice $\Lambda=\mathbb{Z}^n$ be equipped with a skew-symmetric integer form $\langle \cdot, \cdot \rangle$. Introduce the $q$-multiplication operation in the vector space $\Upsilon=\operatorname{Span}\{ Z_\lambda\}_{\lambda\in\Lambda}$ as follows: $Z_\lambda Z_\mu=q^{\langle \lambda,\mu\rangle} Z_{\lambda+\mu}$.  The algebra $\Upsilon$ is called a \emph{quantum torus}. Fixing a basis $\{e_i\}$ in $\Lambda$, we consider $\Upsilon$ as a non-commutative algebra of Laurent polynomials in variables $Z_i:=Z_{e_i}$, $i\in [1,N]$. For any sequence ${\bf a}=(a_1,\dots,a_t)$, $a_i\in [1,N]$, let $\Pi_{\bf a}$ denote the monomial $\Pi_{\bf a}=Z_{a_1} Z_{a_2} \dots Z_{a_t}$. Let $\lambda_{\bf a}=\sum_{j=1}^t e_{a_j}$. Element $Z_{\lambda_{\bf a}}$ is called in physical literature the Weyl form of $\Pi_{\bf a}$ and we denote it by two-sided colons $\col{\Pi_{\bf a}}$  It is easy to see that  $\col{\Pi_{\bf a}}=Z_{\lambda_{\bf a}}=q^{-\sum_{j<k}\langle e_{a_j}, e_{a_k}\rangle} \Pi_{\bf a}$. 

Every $k\times p$ (quantum) matrix $M$ lies in the direct product of $\hbox{Mat}_{k\times p}\otimes \bf Q$, where $\bf Q$ is the quantum operatorial space shared by all matrices and all matrix elements. We can interpret the matrix $M$ as a morphism from ${\mathbb C}^p \otimes W$ into a vector space ${\mathbb C}^k\otimes W$, where $W$ is an irreducible representation space for  operators from the set $\bf Q$.  In what follows, we assume that the order of product of quantum operators in $\bf Q$ always coincides with the order of the product of  corresponding matrix elements.  All {\em classical} operators act  as the unit operators in $W$; in particular, all $R$-matrices are assumed to be classical. Since we assume the representation $W$ to be irreducible, all Casimirs are classical operators. To shorten the writing, we usually omit direct product symbols instead indicating by the index above the symbol of an operator the number of the space ${\mathbb C}^k$ in which this operator acts nontrivially. For example, $\sheet{1} M\sheet{2} V$ of two operators  $M=\sum_{i,j} e_{ij}\otimes m_{ij}$ and $V=\sum_{k,l}e_{kl}\otimes v_{kl}$ with $m_{ij}$ and $v_{kl}$ from the set $\bf Q$ denotes the operator acting in ${\mathbb C}^k\otimes {\mathbb C}^k \otimes W$
$$
\sheet{1}{M}\sheet{2} V :=\sum_{i,j,k,l}\bigl(\sheet{1}{e_{ij}}\otimes \sheet{2} e_{kl}\otimes m_{ij}v_{kl} \bigr),
$$
where we use the standard notation $e_{ij}$ for an elementary matrix whose all elements vanish except the unit element at the intersection of the $i$th row and $j$th column. 

In what follows, we let $I$ denote the unit matrix acting on ${\mathbb C}^k$, $\mathbb I$ the unit operator acting in $W$, and $\mathbf I=I\otimes\mathbb I$ the unit operator acting in the direct product ${\mathbb C}^k\otimes W$. We imply that $\sheet{1}{M}:=\sheet{1}{M}\otimes \sheet{2}{\mathbf I}$. The operation of matrix inversion, $M^{-1}$, acts in the both spaces: classical and quantum, so $M^{-1}M=\mathbf I$ and the order of operators in the quantum space $\bf Q$ is inverted: e.g., $\bigl[\sheet{1}{M}\sheet{2} V\bigr]^{-1}=\sheet{2} V{}^{-1} \sheet{1} M{}^{-1}$; on the contrary, transposition acts only in the classical space, so $\bigl[\sheet{1}{M}\sheet{2} V\bigr]^{\text{T}}=\sheet{1} M{}^{\text{T}} \sheet{2} V{}^{\text{T}}$. Finally, the Hermitian conjugation again permutes operators,  $\bigl[\sheet{1}{M}\sheet{2} V\bigr]^{\dagger}=\sheet{1} M{}^{\dagger} \sheet{2} V{}^{\dagger}$.


\subsection{Quantum transport matrices and Fock-Goncharov coordinates}\label{sec:QFG} 
We describe now how quantized transport matrices are expressed in terms of quantized Fock--Goncharov parameters.

\begin{figure}[h]
	\begin{pspicture}(-4,-2.5)(4,4){
		\newcommand{\LEFTDOWNARROW}{%
			{\psset{unit=1}
				\rput(0,0){\psline[doubleline=true,linewidth=1pt, doublesep=1pt, linecolor=black]{<-}(0,0)(.765,.45)}
		}}
		\newcommand{\DOWNARROW}{%
	{\psset{unit=1}
					\rput(0,0){\psline[doubleline=true,linewidth=1pt, doublesep=1pt, linecolor=black]{->}(0,0.1)(0,-0.566)}
		\put(0,0){\pscircle[fillstyle=solid,fillcolor=white]{.1}}
}}
		\newcommand{\LEFTUPARROW}{%
	{\psset{unit=1}
		\rput(0,0){\psline[doubleline=true,linewidth=1pt, doublesep=1pt, linecolor=black]{->}(0,0)(-.765,.45)}
}}
	\newcommand{\STARUP}{
			{\psset{unit=1}
	\rput(0,0){\psline[doubleline=true,linewidth=1pt, doublesep=1pt, linecolor=black]{<-}(0.06,-0.03)(.4,-.24)}
	\rput(0,0){\psline[doubleline=true,linewidth=1pt, doublesep=1pt, linecolor=black]{<-}(0,0.1)(0,.466)}
	\rput(0,0){\psline[doubleline=true,linewidth=1pt, doublesep=1pt, linecolor=black]{->}(0,0)(-.45,-.23)}
	\put(0,0){\pscircle[fillstyle=solid,fillcolor=black]{.1}}
	\put(0,.566){\pscircle[fillstyle=solid,fillcolor=white]{.1}}
}}
		\newcommand{\PATGEN}{%
			{\psset{unit=1}
				\rput(0,0){\psline[linecolor=blue,linewidth=2pt]{->}(0,0)(.45,.765)}
				\rput(0,0){\psline[linecolor=blue,linewidth=2pt]{->}(1,0)(0.1,0)}
				\rput(0,0){\psline[linecolor=blue,linewidth=2pt]{->}(0,0)(.45,-.765)}
				\put(0,0){\pscircle[fillstyle=solid,fillcolor=red]{.1}}
		}}
		\newcommand{\PATLEFT}{%
			{\psset{unit=1}
				\rput(0,0){\psline[linecolor=blue,linewidth=2pt,linestyle=dashed]{->}(0,0)(.45,.765)}
				\rput(0,0){\psline[linecolor=blue,linewidth=2pt]{->}(1,0)(0.1,0)}
				\rput(0,0){\psline[linecolor=blue,linewidth=2pt]{->}(0,0)(.45,-.765)}
				\put(0,0){\pscircle[fillstyle=solid,fillcolor=red]{.1}}
		}}
		\newcommand{\PATRIGHT}{%
			{\psset{unit=1}
				\rput(0,0){\psline[linecolor=blue,linewidth=2pt,linestyle=dashed]{->}(0,0)(.45,-.765)}
				\put(0,0){\pscircle[fillstyle=solid,fillcolor=lightgray]{.1}}
		}}
		\newcommand{\PATBOTTOM}{%
			{\psset{unit=1}
				\rput(0,0){\psline[linecolor=blue,linewidth=2pt]{->}(0,0)(.45,.765)}
				\rput(0,0){\psline[linecolor=blue,linewidth=2pt,linestyle=dashed]{->}(1,0)(0.1,0)}
				\put(0,0){\pscircle[fillstyle=solid,fillcolor=red]{.1}}
		}}
		\newcommand{\PATTOP}{%
			{\psset{unit=1}
				\rput(0,0){\psline[linecolor=blue,linewidth=2pt]{->}(1,0)(0.1,0)}
				\rput(0,0){\psline[linecolor=blue,linewidth=2pt]{->}(0,0)(.45,-.765)}
				\put(0,0){\pscircle[fillstyle=solid,fillcolor=red]{.1}}
		}}
		\newcommand{\PATBOTRIGHT}{%
			{\psset{unit=1}
				\rput(0,0){\psline[linecolor=blue,linewidth=2pt]{->}(0,0)(.45,.765)}
				\put(0,0){\pscircle[fillstyle=solid,fillcolor=red]{.1}}
				\put(.5,0.85){\pscircle[fillstyle=solid,fillcolor=lightgray]{.1}}
		}}
		\newcommand{\PATNORTH}{%
			{\psset{unit=1}
				\rput(0,0){\psline[linecolor=blue,linewidth=2pt,linestyle=dashed]{->}(-.45,-.765)(0,0)}
				\rput(0,0){\psline[linecolor=blue,linewidth=2pt,linestyle=dashed]{<-}(.45,-.765)(0,0)}
				\put(0,0){\pscircle[fillstyle=solid,fillcolor=lightgray]{.1}}
		}}
		\newcommand{\PATSW}{%
			{\psset{unit=1}
				\rput(0,0){\psline[linecolor=blue,linewidth=2pt,linestyle=dashed]{->}(0,0)(.45,.765)}
				\rput(0,0){\psline[linecolor=blue,linewidth=2pt,linestyle=dashed]{<-}(0.1,0)(0.9,0)}
				\put(0,0){\pscircle[fillstyle=solid,fillcolor=lightgray]{.1}}
		}}
		\newcommand{\PATSE}{%
			{\psset{unit=1}
				\rput(0,0){\psline[linecolor=blue,linewidth=2pt,linestyle=dashed]{<-}(0,0)(-.45,.765)}
				\rput(0,0){\psline[linecolor=blue,linewidth=2pt,linestyle=dashed]{<-}(-0.9,0)(-0.1,0)}
				\put(0,0){\pscircle[fillstyle=solid,fillcolor=lightgray]{.1}}
		}}
		\multiput(-2.5,-0.85)(0.5,0.85){4}{\PATLEFT}
		\multiput(-2,-1.7)(1,0){4}{\PATBOTTOM}
		\put(-0.5,2.55){\PATTOP}
		\put(0,3.4){\PATNORTH}
		\multiput(-1.5,-0.85)(1,0){4}{\PATGEN}
		\multiput(-1,0)(1,0){3}{\PATGEN}
		\multiput(-.5,0.85)(1,0){2}{\PATGEN}
		\put(0,1.7){\PATGEN}
		\multiput(-1.5,-0.85)(1,0){4}{\PATGEN}
		\multiput(0.5,2.55)(0.5,-0.85){4}{\PATRIGHT}
		\put(2,-1.7){\PATBOTRIGHT}
		\put(-3,-1.7){\PATSW}
		\put(3,-1.7){\PATSE}
		\multiput(-2,-1.176)(1.0,0){5}{\STARUP}
		\multiput(-1.5,-0.335)(1.0,0){4}{\STARUP}
		\multiput(-1.0,0.5)(1.0,0){3}{\STARUP}
		\multiput(-.5,1.4)(1.0,0){2}{\STARUP}
		\put(0,2.3){\STARUP}
		\multiput(2.6,-1.4)(-0.5,.85){6}{\LEFTDOWNARROW}
		\multiput(-2.6,-1.4)(0.5,.85){6}{\LEFTUPARROW}
		\multiput(-2.5,-1.5)(1.0,0){6}{\DOWNARROW}
		\multiput(-0.5,2.55)(0.5,-0.85){6}{\pscircle[fillstyle=solid,fillcolor=violet]{.1}}
		\multiput(-1,1.7)(0.5,-0.85){5}{\pscircle[fillstyle=solid,fillcolor=purple]{.1}}
		\multiput(-1.5,0.85)(0.5,-0.85){4}{\pscircle[fillstyle=solid,fillcolor=magenta]{.1}}
		\multiput(-2,0)(0.5,-0.85){3}{\pscircle[fillstyle=solid,fillcolor=red]{.1}}
		\multiput(-2.5,-0.85)(0.5,-0.85){2}{\pscircle[fillstyle=solid,fillcolor=orange]{.1}}
		\put(1.2,3.2){\makebox(0,0)[br]{\hbox{{$1$}}}}
		\put(1.7,2.4){\makebox(0,0)[br]{\hbox{{$2$}}}}
		\put(2.2,1.6){\makebox(0,0)[br]{\hbox{{$3$}}}}
		\put(2.7,0.8){\makebox(0,0)[br]{\hbox{{$4$}}}}
		\put(3.2,-0.1){\makebox(0,0)[br]{\hbox{{$5$}}}}
		\put(3.7,-1.0){\makebox(0,0)[br]{\hbox{{$6$}}}}

		\put(-0.9,3.2){\makebox(0,0)[br]{\hbox{{$1'$}}}}
		\put(-1.4,2.4){\makebox(0,0)[br]{\hbox{{$2'$}}}}
		\put(-1.9,1.6){\makebox(0,0)[br]{\hbox{{$3'$}}}}
		\put(-2.4,0.8){\makebox(0,0)[br]{\hbox{{$4'$}}}}
		\put(-2.9,-0.1){\makebox(0,0)[br]{\hbox{{$5'$}}}}
		\put(-3.4,-1.0){\makebox(0,0)[br]{\hbox{{$6'$}}}}
		
		\put(-2.2,-2.4){\makebox(0,0)[br]{\hbox{{$1''$}}}}
		\put(-1.2,-2.4){\makebox(0,0)[br]{\hbox{{$2''$}}}}
		\put(-0.2,-2.4){\makebox(0,0)[br]{\hbox{{$3''$}}}}
		\put(0.8,-2.4){\makebox(0,0)[br]{\hbox{{$4''$}}}}
		\put(1.8,-2.4){\makebox(0,0)[br]{\hbox{{$5''$}}}}
		\put(2.8,-2.4){\makebox(0,0)[br]{\hbox{{$6''$}}}}
\put(-2.7,-1.9){\makebox(0,0)[tr]{\hbox{{\small $Z_{(6,0,0)}$}}}}
\put(2.9,-1.9){\makebox(0,0)[tl]{\hbox{{\small $Z_{(0,6,0)}$}}}}
\put(0,3.6){\makebox(0,0)[bc]{\hbox{{\small $Z_{(0,0,6)}$}}}}
	}
	\end{pspicture}
	\begin{pspicture}(-4,-2.5)(4,4){
		\newcommand{\LEFTDOWNARROW}{%
			{\psset{unit=1}
				\rput(0,0){\psline[doubleline=true,linewidth=1pt, doublesep=1pt, linecolor=black]{->}(0,0)(.765,.45)}
		}}
		\newcommand{\DOWNARROW}{%
	{\psset{unit=1}
					\rput(0,0){\psline[doubleline=true,linewidth=1pt, doublesep=1pt, linecolor=black]{->}(0,0.1)(0,-0.566)}
		\put(0,0){\pscircle[fillstyle=solid,fillcolor=white]{.1}}
}}
		\newcommand{\LEFTUPARROW}{%
	{\psset{unit=1}
		\rput(0,0){\psline[doubleline=true,linewidth=1pt, doublesep=1pt, linecolor=black]{<-}(0,0)(-.765,.45)}
}}
	\newcommand{\STARUP}{
			{\psset{unit=1}
	\rput(0,0){\psline[doubleline=true,linewidth=1pt, doublesep=1pt, linecolor=black]{->}(0.06,-0.03)(.4,-.24)}
	\rput(0,0){\psline[doubleline=true,linewidth=1pt, doublesep=1pt, linecolor=black]{<-}(0,0.1)(0,.466)}
	\rput(0,0){\psline[doubleline=true,linewidth=1pt, doublesep=1pt, linecolor=black]{<-}(0,0)(-.45,-.23)}
	\put(0,0){\pscircle[fillstyle=solid,fillcolor=black]{.1}}
	\put(0,.566){\pscircle[fillstyle=solid,fillcolor=white]{.1}}
}}
		\newcommand{\PATGEN}{%
			{\psset{unit=1}
				\rput(0,0){\psline[linecolor=blue,linewidth=2pt]{->}(0,0)(.45,.765)}
				\rput(0,0){\psline[linecolor=blue,linewidth=2pt]{->}(1,0)(0.1,0)}
				\rput(0,0){\psline[linecolor=blue,linewidth=2pt]{->}(0,0)(.45,-.765)}
				\put(0,0){\pscircle[fillstyle=solid,fillcolor=red]{.1}}
		}}
		\newcommand{\PATLEFT}{%
			{\psset{unit=1}
				\rput(0,0){\psline[linecolor=blue,linewidth=2pt,linestyle=dashed]{->}(0,0)(.45,.765)}
				\rput(0,0){\psline[linecolor=blue,linewidth=2pt]{->}(1,0)(0.1,0)}
				\rput(0,0){\psline[linecolor=blue,linewidth=2pt]{->}(0,0)(.45,-.765)}
				\put(0,0){\pscircle[fillstyle=solid,fillcolor=red]{.1}}
		}}
		\newcommand{\PATRIGHT}{%
			{\psset{unit=1}
				\rput(0,0){\psline[linecolor=blue,linewidth=2pt,linestyle=dashed]{->}(0,0)(.45,-.765)}
				\put(0,0){\pscircle[fillstyle=solid,fillcolor=lightgray]{.1}}
		}}
		\newcommand{\PATBOTTOM}{%
			{\psset{unit=1}
				\rput(0,0){\psline[linecolor=blue,linewidth=2pt]{->}(0,0)(.45,.765)}
				\rput(0,0){\psline[linecolor=blue,linewidth=2pt,linestyle=dashed]{->}(1,0)(0.1,0)}
				\put(0,0){\pscircle[fillstyle=solid,fillcolor=red]{.1}}
		}}
		\newcommand{\PATTOP}{%
			{\psset{unit=1}
				\rput(0,0){\psline[linecolor=blue,linewidth=2pt]{->}(1,0)(0.1,0)}
				\rput(0,0){\psline[linecolor=blue,linewidth=2pt]{->}(0,0)(.45,-.765)}
				\put(0,0){\pscircle[fillstyle=solid,fillcolor=red]{.1}}
		}}
		\newcommand{\PATBOTRIGHT}{%
			{\psset{unit=1}
				\rput(0,0){\psline[linecolor=blue,linewidth=2pt]{->}(0,0)(.45,.765)}
				\put(0,0){\pscircle[fillstyle=solid,fillcolor=red]{.1}}
				\put(.5,0.85){\pscircle[fillstyle=solid,fillcolor=lightgray]{.1}}
		}}
		\newcommand{\PATNORTH}{%
			{\psset{unit=1}
				\rput(0,0){\psline[linecolor=blue,linewidth=2pt,linestyle=dashed]{->}(-.45,-.765)(0,0)}
				\rput(0,0){\psline[linecolor=blue,linewidth=2pt,linestyle=dashed]{<-}(.45,-.765)(0,0)}
				\put(0,0){\pscircle[fillstyle=solid,fillcolor=lightgray]{.1}}
		}}
		\newcommand{\PATSW}{%
			{\psset{unit=1}
				\rput(0,0){\psline[linecolor=blue,linewidth=2pt,linestyle=dashed]{->}(0,0)(.45,.765)}
				\rput(0,0){\psline[linecolor=blue,linewidth=2pt,linestyle=dashed]{<-}(0.1,0)(0.9,0)}
				\put(0,0){\pscircle[fillstyle=solid,fillcolor=lightgray]{.1}}
		}}
		\newcommand{\PATSE}{%
			{\psset{unit=1}
				\rput(0,0){\psline[linecolor=blue,linewidth=2pt,linestyle=dashed]{<-}(0,0)(-.45,.765)}
				\rput(0,0){\psline[linecolor=blue,linewidth=2pt,linestyle=dashed]{<-}(-0.9,0)(-0.1,0)}
				\put(0,0){\pscircle[fillstyle=solid,fillcolor=lightgray]{.1}}
		}}
		\multiput(-2.5,-0.85)(0.5,0.85){4}{\PATLEFT}
		\multiput(-2,-1.7)(1,0){4}{\PATBOTTOM}
		\put(-0.5,2.55){\PATTOP}
		\put(0,3.4){\PATNORTH}
		\multiput(-1.5,-0.85)(1,0){4}{\PATGEN}
		\multiput(-1,0)(1,0){3}{\PATGEN}
		\multiput(-.5,0.85)(1,0){2}{\PATGEN}
		\put(0,1.7){\PATGEN}
		\multiput(-1.5,-0.85)(1,0){4}{\PATGEN}
		\multiput(0.5,2.55)(0.5,-0.85){4}{\PATRIGHT}
		\put(2,-1.7){\PATBOTRIGHT}
		\put(-3,-1.7){\PATSW}
		\put(3,-1.7){\PATSE}
		\multiput(-2,-1.176)(1.0,0){5}{\STARUP}
		\multiput(-1.5,-0.335)(1.0,0){4}{\STARUP}
		\multiput(-1.0,0.5)(1.0,0){3}{\STARUP}
		\multiput(-.5,1.4)(1.0,0){2}{\STARUP}
		\put(0,2.3){\STARUP}
		\multiput(2.6,-1.4)(-0.5,.85){6}{\LEFTDOWNARROW}
		\multiput(-2.6,-1.4)(0.5,.85){6}{\LEFTUPARROW}
		\multiput(-2.5,-1.5)(1.0,0){6}{\DOWNARROW}
		\multiput(-0.5,2.55)(0.5,-0.85){6}{\pscircle[fillstyle=solid,fillcolor=violet]{.1}}
		\multiput(-1,1.7)(0.5,-0.85){5}{\pscircle[fillstyle=solid,fillcolor=purple]{.1}}
		\multiput(-1.5,0.85)(0.5,-0.85){4}{\pscircle[fillstyle=solid,fillcolor=magenta]{.1}}
		\multiput(-2,0)(0.5,-0.85){3}{\pscircle[fillstyle=solid,fillcolor=red]{.1}}
		\multiput(-2.5,-0.85)(0.5,-0.85){2}{\pscircle[fillstyle=solid,fillcolor=orange]{.1}}
		\put(1.2,3.2){\makebox(0,0)[br]{\hbox{{$1$}}}}
		\put(1.7,2.4){\makebox(0,0)[br]{\hbox{{$2$}}}}
		\put(2.2,1.6){\makebox(0,0)[br]{\hbox{{$3$}}}}
		\put(2.7,0.8){\makebox(0,0)[br]{\hbox{{$4$}}}}
		\put(3.2,-0.1){\makebox(0,0)[br]{\hbox{{$5$}}}}
		\put(3.7,-1.0){\makebox(0,0)[br]{\hbox{{$6$}}}}

		\put(-0.9,3.2){\makebox(0,0)[br]{\hbox{{$1'$}}}}
		\put(-1.4,2.4){\makebox(0,0)[br]{\hbox{{$2'$}}}}
		\put(-1.9,1.6){\makebox(0,0)[br]{\hbox{{$3'$}}}}
		\put(-2.4,0.8){\makebox(0,0)[br]{\hbox{{$4'$}}}}
		\put(-2.9,-0.1){\makebox(0,0)[br]{\hbox{{$5'$}}}}
		\put(-3.4,-1.0){\makebox(0,0)[br]{\hbox{{$6'$}}}}

		\put(-2.2,-2.4){\makebox(0,0)[br]{\hbox{{$1''$}}}}
		\put(-1.2,-2.4){\makebox(0,0)[br]{\hbox{{$2''$}}}}
		\put(-0.2,-2.4){\makebox(0,0)[br]{\hbox{{$3''$}}}}
		\put(0.8,-2.4){\makebox(0,0)[br]{\hbox{{$4''$}}}}
		\put(1.8,-2.4){\makebox(0,0)[br]{\hbox{{$5''$}}}}
		\put(2.8,-2.4){\makebox(0,0)[br]{\hbox{{$6''$}}}}
\put(-2.7,-1.9){\makebox(0,0)[tr]{\hbox{{\small $Z_{(6,0,0)}$}}}}
\put(2.9,-1.9){\makebox(0,0)[tl]{\hbox{{\small $Z_{(0,6,0)}$}}}}
\put(0,3.6){\makebox(0,0)[bc]{\hbox{{\small $Z_{(0,0,6)}$}}}}
	}
	\end{pspicture}

	\caption{\small
	Two directed networks $N$ dual to the Fock--Goncharov $b_6$-quiver. Double arrows are edges of the directed network,  cluster variables correspond to faces of $N$ (vertices of the $b_n$-quiver dual to $N$). On the left-hand side we present the directed network defining the transport matrices $M_1$ (from the side $\{1-6\}$ to the side $\{1'-6'\}$ and $M_2$ (from the side $\{1-6\}$ to the side $\{1''-6''\}$) and on the right-hand side we present the directed network defining the transport matrix $M_3$ (from the side $\{1'-6'\}$ to the side $\{1''-6''\})$.  We also indicate three barycentrically enumerated cluster variables at three corners of the quiver.
	}
	\label{fi:triangle}
\end{figure}
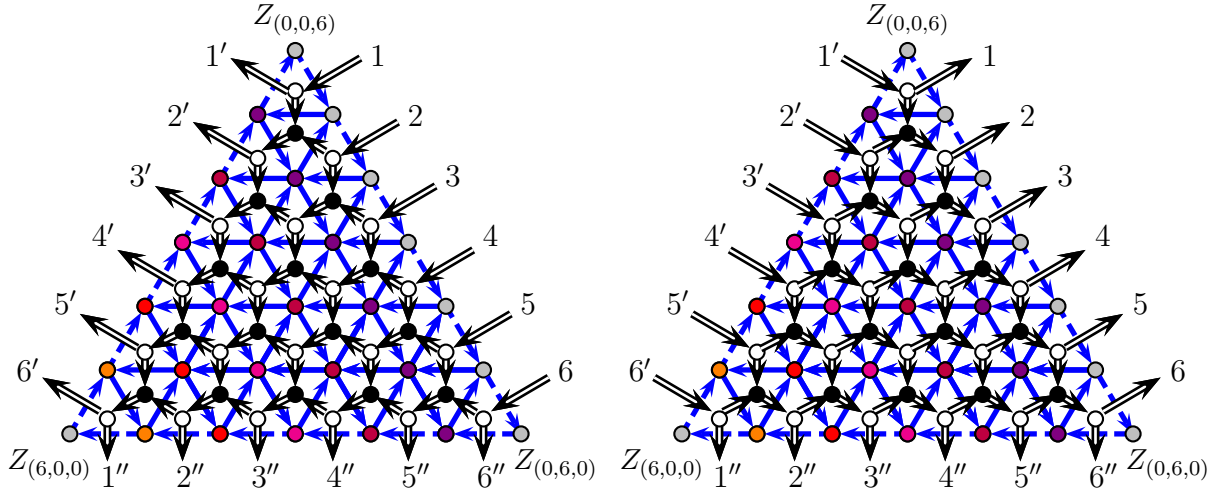

In the quantization of $\X_{SL_n,\Sigma}$ for $\Sigma_{0,1,3}$,  the quantized Fock--Goncharov variables form a quantum torus $\Upsilon$ with commutation relation described by the $b_n$-quiver ($b_n$ stands for a Borel subalgebra of $sl_n$) shown in Fig.~\ref{fi:triangle} by solid and dashed lines. Vertices of the quiver label quantum Fock-Goncharov coordinates $Z_\alpha$  while arrows  encode commutation relations:  if there are $m$ arrows from vertex $\alpha$ to $\beta$ then $Z_{\beta}Z_{\alpha}=q^{-2m} Z_{\alpha}Z_{\beta}$. Dashed arrow counts as $m=1/2$. In particular, a solid arrow from $Z_{\alpha}$ to $Z_{\beta}$ implies $Z_{\beta}Z_{\alpha}=q^{-2} Z_{\alpha}Z_{\beta}$,  a dashed arrow from $Z_{\alpha}$ to $Z_{\beta}$ implies $Z_{\beta}Z_{\alpha}=q^{-1}Z_{\alpha}Z_{\beta}$, and, say, a double arrow from $Z_\alpha$ to $Z_\beta$ means $Z_{\beta}Z_{\alpha}=q^{-4}Z_{\alpha}Z_{\beta}$.  Vertices not connected by an arrow commute.

We next construct a {\sl directed network} $N$ dual to the $b_n$-quiver. 
Dual directed network is a directed embedded graph in the disk (triangle) whose \emph{three-valent} vertices are colored black and white and orientation is compatible with the coloring, i.e. each white vertex has exactly one incoming edge while every black vertex has exactly one outgoing edge. Vertices of the quiver that is dual to a directed network are located in the faces of the directed network (one vertex per each face) while directed edges of the quiver intersect black-white edges of the network in such a way that quiver's arrows travel clockwise around white vertices of network and counterclockwise around black vertices.  Since the quiver orientation depends only on a pattern of black and white vertices and not on orientations of network edges, several directed networks may correspond to the same quiver, but in the considered case orientation of the dual network is determined by fixing $n$ sources and $2n$ sinks among boundary vertices. In particular, choosing all sources of the network located along one side of the triangle Fig.~\ref{fi:triangle} we fix the dual directed network uniquely.


We therefore have several possible choices of directed coloring-compatible networks for the $b_n$-quiver, two of which relevant for our studies are depicted in the figure: every solid arrow of a quiver has a white vertex on the right and a black vertex on the left;  all dashed arrows have white vertices on the right. Directed edges of the network (double lines in the figure) are dual to edges of the quiver with sources and sinks dual to dashed arrows; it is easy to see that the choice of sources and sinks then determines the pattern of directed double lines inside the graph in a unique way.  

We adopt a barycentric enumeration of vertices and cluster variables $Z_{(i,j,k)}$ of the $b_n$-quiver: $i,j,k$ are nonnegative integers such that $i+j+k=n$ (one can think about these vertices as integer points of intersection of the lattice $\mathbb Z_{+,0}\otimes \mathbb Z_{+,0}\otimes \mathbb Z_{+,0}$ with the plane $x+y+z=n$). 
 
We assign to every oriented path $P:j\leadsto i'$ from the right side to the left side or to the bottom side $P:j\leadsto i''$ the {\em quantum weight}
\be\label{eq:pathweight}
w(P)=\col{\mathop{\prod_{\text{ faces }\alpha\text{ lie to the right}}}_{\text{ of the path } P} Z_\alpha}.
\ee

\begin{definition}\label{df:M}
We define two $n\times n$ \emph{non-normalized quantum transition matrices} 
$$
(\mathcal M_1)_{i,j}=\mathop{\sum_{\text{directed path\,}{P: j\leadsto i'}}}_{\text{ from right to left}} w(P)\quad \hbox{and}\quad
(\mathcal M_2)_{i,j}=\mathop{\sum_{\text{directed path\,}{P: j\leadsto i''}}}_{\text{ from right to bottom}} w(P).
$$
Note that ${\mathcal M}_1$ is a lower-triangular matrix and ${\mathcal M}_2$ is an upper-triangular matrix.
\end{definition}

Theorem 9.3 of \cite{ChSh2} states that for any planar network with separated $n$ sources and $m$ sinks, defining entries of $m\times n$ quantum transition matrix as in Definition~\ref{df:M}, we obtain that these entries enjoy the {\em quantum $R$-matrix relation} 
\be\label{RMM}
\mathcal R_{m}(q) \sheet{1}{\mathcal M}\otimes \sheet{2}{\mathcal M}=\sheet{2}{\mathcal M}\otimes \sheet{1}{\mathcal M}\mathcal R_n(q),
\ee
where $\mathcal R_{k}(q)$ denotes a $k^2\times k^2$ matrix
\be\label{R-matrix1}
\mathcal R_{k}(q)= \sum_{1\le i,j\le k} \sheet{1}e_{ii}\otimes \sheet{2} e_{jj}+(q-1)\sum_{1\le i\le k} \sheet{1}e_{ii}\otimes \sheet{2} e_{ii}
+(q-q^{-1})\sum_{1\le j<i \le k} \sheet{1}e_{ij}\otimes \sheet{2} e_{ji}
\ee
having the properties that
$\mathcal R^{-1}_k(q)=\mathcal R_k(q^{-1})$ and $\mathcal R_k(q)-\mathcal R^{\text{T}}_k(q^{-1})=(q-q^{-1})P_k$,
with $P_k$ the standard permutation matrix $P_k:=\sum_{1\le i,j\le k}\sheet{1}e_{ij}\otimes \sheet{2} e_{ji}$. Note that the total transposition of $\mathcal R_{k}(q)$ results in interchanging the space labels $1\leftrightarrow 2$.

In Fig.~\ref{fi:triangle} we have an example of a directed network with $n=6$ sources
and $2n=12$ sinks; for any $n$ we present the corresponding quantum transition matrix in the block form
$\mathcal M=\begin{pmatrix} \mathcal M_1 \\ \mathcal M_2\end{pmatrix}$ with $\mathcal M_1$ and $\mathcal M_2$ denoting the respective upper and lower $n\times n$ blocks. Relations (\ref{RMM}) then imply the $R$-matrix commutation relations between these blocks: 
\be\label{R1}
\mathcal R_n(q)\sheet{1}{\mathcal M}_i \otimes \sheet{2}{\mathcal M}_i=\sheet{2}{\mathcal M}_i\otimes \sheet{1}{\mathcal M}_i \mathcal R_{n}(q),\quad i=1,2,
\ee
and
\be\label{R2}
\sheet{1}{\mathcal M}_1 \otimes \sheet{2}{\mathcal M}_2=\sheet{2}{\mathcal M}_2\otimes \sheet{1}{\mathcal M}_1 \mathcal R_{n}(q).
\ee
We also introduce the third transition matrix $\mathcal M_3$ for paths going from the side $\{1'-6'\}$ to the side $\{1''-6''\})$ in the transformed network in the right side of Fig.~\ref{fi:triangle}.

We now define the  \emph{quantum transport} matrices of the \emph{standard $b_n$-quiver}:

\begin{definition}\label{def:M}
Quantum transport matrices 
read
$$
M_1=QS{\mathcal M}_1,\quad  M_2=QS{\mathcal M}_2, \text{ and }  M_3=QS{\mathcal M}_3,
$$ 
where  $Q:=\sum_{i=1}^n (q)^{-i+1/2}e_{i,i}\otimes \mathbb I$ is a diagonal matrix and $\displaystyle S=\sum_{i=1}^n (-1)^{i+1}e_{i,n+1-i}\otimes \mathbb I$  is an antidiagonal matrix.
\end{definition}

Note that  
$$
\sheet{1}S\otimes \sheet{2}S \mathcal R_n(q)=\mathcal R^{\text{T}}_n(q)\sheet{1}S\otimes \sheet{2}S\  \hbox{for any antidiagonal classical matrix $S$}. 
$$
and
$$
\sheet{1}Q\otimes \sheet{2}Q \mathcal R_n(q)=\mathcal R_n(q)\sheet{1}Q\otimes \sheet{2}Q\  \hbox{for any diagonal classical matrix $Q$};
$$
in particular, setting $Q=AB$ with two diagonal classical invertible matrices $A$ and $B$, we have
\be\label{ABR}
\sheet{1}A\otimes \sheet{2}B \mathcal R_n(q) \sheet{1}B{}^{-1}\otimes \sheet{2}A{}^{-1}= \sheet{1}B{}^{-1}\otimes \sheet{2}A{}^{-1} \mathcal R_n(q)\sheet{1}A\otimes \sheet{2}B.
\ee

\begin{remark}\label{rem:QS}
Multiplying a transition matrix by the matrix $S$ is standard in the Fock--Goncharov approach in order to  make orientations of boundaries of adjacent triangles compatible and hence to ``invert'' the corresponding flags upon reaching a boundary of a network; additional $q$-factors collected in the matrix $Q$ are necessary for ensuring the groupoid property below; we do not normalize transport matrices (and, correspondingly, the $R$-matrix) as in \cite{ChSh2} because this normalization does not affect our consideration below. We also omit the sign of the direct product  in the subsequent text.
\end{remark}

We have the following theorem.

\begin{theorem}\label{th:MM}\cite{ChSh2}.
The quantum transport matrices $M_1$,  $M_2$, and $M_3$ satisfy the relations
\begin{align*}
&{\mathcal R}^{\text{T}}_{n}(q)\sheet{1}M_i  \sheet{2} M_i =\sheet{2}M_i \sheet{1}M_i {\mathcal  R}_{n}(q),\quad i=1,2,3, \\
&\sheet{1}M_1  \sheet{2} M_2=\sheet{2}M_2 \sheet{1}M_1 {\mathcal R}_{n}(q),\quad \sheet{1}M_1  \sheet{2} M_3=\sheet{2}M_3 {\mathcal R}^{\text{T}}_{n}(q) \sheet{1}M_1,\quad  \sheet{1}M_2  \sheet{2} M_3= {\mathcal R}_{n}(q) \sheet{2}M_3 \sheet{1}M_2 
\end{align*}
with the quantum trigonometric $R$-matrix (\ref{R-matrix1}).
\end{theorem}

\begin{theorem}\label{th:groupoid}\cite{ChSh2}.
The quantum transport matrices satisfy the {\em quantum groupoid condition}
$$M_3 M_1= M_2.$$
This condition is consistent with the quantum commutation relations in Theorem~\ref{th:MM}.
\end{theorem}

Using relations (\ref{ABR}) we obtain
\begin{lemma}\label{lm:symmetries}
 The commutation relations in Theorem~\ref{th:MM}, as well as the groupoid condition in Theorem~\ref{th:groupoid} are invariant under transformations
$$
M_1\to AM_1C,\quad M_2\to BM_2C, \quad M_3\to B M_3 A^{-1},
$$ 
where $A$, $B$, and $C$ are arbitrary nondegenerate classical diagonal matrices.
\end{lemma}

\section{The groupoid of upper triangular matrices}\label{s:A}
\setcounter{equation}{0}

\subsection{Representing an upper-triangular $\mathbb A$}

Consider a special combination of $M_1$ and $M_2$:
\be\label{A}
\mathbb A:=M_1^{\text{T}}D M_2=M_1^{\text{T}} D M_3 M_1,
\ee
where $D$ is any classical diagonal matrix. (In notations of Lemma~\ref{lm:symmetries}, $D=A^{\text{T}}B$.) We can use the freedom in choice of $D$ to attain the canonical form (\ref{A-quantum}) of the quantum $\mathbb A$-matrix with the diagonal entries equal to $q^{-1/2}$. Then the properly normalized expression for this matrix reads
\be
\label{eq:A-norm}
\mathbb A^\hbar =\mathcal M_1^{\text{T}}\mathcal M_3 QS \mathcal M_1,
\ee
i.e., $D^{-1}=S^{\text{T}}Q^{\text{T}}QS$,
and we use expression (\ref{eq:A-norm}) when solving the problem of singular values of $\lambda$ in Sec.~\ref{s:character}.

Note that since $M_1$ and $M_1^{\text{T}}$ are upper-anti-diagonal matrices and $M_2$ is a lower-anti-diagonal matrix in the case of $\Sigma_{0,1,3}$, the matrix $\mathbb A$ is  automatically upper-triangular.

The following theorem was proved in \cite{ChSh2} for $\mathbb A$ having the form $M_1^{\text{T}} M_2$. We slightly modify this proof to adjust it to $\mathbb A$ of the form $M_1^{\text{T}} D M_3 M_1$.

\begin{theorem}\label{th:A}
For the matrices $M_1$ and $M_3$ enjoying commutation relations in Theorem~\ref{th:MM} and for any diagonal matrix $D$,
the matrix $\mathbb A$ given by (\ref{A}) enjoys the quantum reflection equation
$$
{\mathcal R}_{n}(q)\sheet{1}{\mathbb A} {\mathcal R}_{n}^{\text{t}_1}(q)\sheet{2} {\mathbb A}=
\sheet{2}{\mathbb A} {\mathcal R}_{n}^{\text{t}_1}(q)\sheet{1} {\mathbb A}{\mathcal R}_{n}(q)
$$
with the trigonometric $R$-matrix (\ref{R-matrix1}), where $R_{n}^{\text{t}_1}(q)$ is a partially transposed (w.r.t. the first space) $R$-matrix.
\end{theorem}

The proof uses only $R$-matrix relations. Transposing relations in Theorem~\ref{th:MM}  with respect to different spaces, we obtain
$$
\sheet{1}M{}_i^{\text{T}} {\mathcal R}_{n}^{\text{t}_2}(q)  \sheet{2} M_i=\sheet{2}M_i  {\mathcal R}_{n}^{\text{t}_1}(q) \sheet{1} M{}_i^{\text{T}} , \quad
{\mathcal R}_{n}(q) \sheet{1}M{}_1^{\text{T}}  \sheet{2}M{}_1^{\text{T}}= \sheet{2}M{}_1^{\text{T}}  \sheet{1}M{}_1^{\text{T}}{\mathcal R}_{n}^{\text{T}} (q),\ \hbox{and}\ 
\sheet{1}M{}_1^{\text{T}} \sheet{2}M{}_3=  \sheet{2}M{}_3  \sheet{1}M{}_1^{\text{T}} {\mathcal R}_{n}^{\text{t}_2}(q).
$$
Note that $ {\mathcal R}_{n}^{\text{t}_1}(q)$ retains its form if we interchange the space indices $1\leftrightarrow 2$; note also that $D$ is a classical matrix, so, say, $\sheet{1}D$ commutes with all $\sheet{2}M_i$. We have
\begin{align}
&{\mathcal R}_{n}(q) \sheet{1}M{}_1^{\text{T}} \sheet{1}D  \sheet{1}M_3\tcr{ \sheet{1}M{}_1 {\mathcal R}_{n}^{\text{t}_1}(q) \sheet{2}M{}_1^{\text{T}} } \sheet{2}D  \sheet{2}M_3 \sheet{2}M_1=
{\mathcal R}_{n}(q) \sheet{1}M{}_1^{\text{T}}  \sheet{1}D \tcr{  \sheet{1}M{}_3 \sheet{2}M{}_1^{\text{T}}  {\mathcal R}_{n}^{\text{t}_2}(q)}  \sheet{2}D \tcb{  \sheet{1}M_1 \sheet{2}M_3} \sheet{2}M_1\nonumber\\
=&\tcr{ {\mathcal R}_{n}(q) \sheet{1}M{}_1^{\text{T}}  \sheet{2}M{}_1^{\text{T}} }  \sheet{1}D  \sheet{2}D \sheet{1}M{}_3 \sheet{2}M{}_3 \tcb{ {\mathcal R}_{n}^{\text{T}} (q)  \sheet{1}M{}_1  \sheet{2}M{}_1}
=  \sheet{2}M{}_1^{\text{T}}  \sheet{1}M{}_1^{\text{T}} \tcr{ {\mathcal R}_{n}^{\text{T}} (q)  \sheet{1}D  \sheet{2}D }  \sheet{1}M{}_3  \sheet{2}M{}_3  \sheet{2}M{}_1  \sheet{1}M{}_1 {\mathcal R}_{n}(q) \nonumber\\
=& \sheet{2}M{}_1^{\text{T}}  \sheet{1}M{}_1^{\text{T}}  \sheet{1}D  \sheet{2}D \tcr{ {\mathcal R}_{n}^{\text{T}} (q)   \sheet{1}M{}_3  \sheet{2}M{}_3 }  \sheet{2}M{}_1  \sheet{1}M{}_1 {\mathcal R}_{n}(q) 
=\sheet{2}M{}_1^{\text{T}}   \sheet{2}D \tcr{  \sheet{1}M{}_1^{\text{T}}  \sheet{2}M{}_3} \sheet{1}D  \tcb{ \sheet{1}M{}_3  {\mathcal R}_{n}(q) \sheet{2}M{}_1} \sheet{1}M{}_1 {\mathcal R}_{n}(q)\label{long}\\
=& \sheet{2}M{}_1^{\text{T}}  \sheet{2}D  \sheet{2}M{}_3 \tcr{  \sheet{1}M{}_1^{\text{T}} {\mathcal R}_{n}^{\text{t}_2}(q) \sheet{2}M{}_1}  \sheet{1}D \sheet{1}M{}_3 \sheet{1}M{}_1 {\mathcal R}_{n}(q)
= \sheet{2}M{}_1^{\text{T}}  \sheet{2}D  \sheet{2}M_3 \sheet{2}M{}_1 {\mathcal R}_{n}^{\text{t}_1}(q) \sheet{1}M{}_1^{\text{T}}  \sheet{1}D  \sheet{1}M_3 \sheet{1}M_1 {\mathcal R}_{n}(q),\nonumber
\end{align}
which completes the proof.

We have therefore {\em a Darboux coordinate representation} for operators satisfying the reflection equation. Moreover, all matrix elements of $\mathbb A$ are Laurent polynomials with positive coefficients of $Z_\alpha$ and $q$. 

By construction of quantum transport matrices in Sec.~\ref{sec:QFG}, all matrix elements of $M_1$ and $M_2$ are Weyl-ordered. For $[\mathbb A]_{i,j}=\sum\limits_{k=i}^j [M_1]_{k,i}[M_2]_{k,j}$ we obtain that for $i<j$, $[M_1]_{k,i}$ commutes with  $[M_2]_{k,j}$ (all paths contributing to $[M_1]_{k,i}$ are disjoint from all paths contributing to $[M_2]_{k,j}$ for $i<j$), so the corresponding products are automatically Weyl-ordered,  $[\mathbb A]_{i,j}=\sum\limits_{k=i}^j \col{[M_1]_{k,i}[M_2]_{k,j}}$. For $i=j$, the corresponding two paths share the common starting edge, and then $[\mathbb A]_{i,i}=q^{-1/2} \col{[M_1]_{i,i}[M_2]_{i,i}}$. This explains the appearance of $q^{-1/2}$ factors on the diagonal of the quantum matrix $\mathbb A^\hbar$ (see (\ref{A-quantum}), \cite{ChM}). Note that all Weyl-ordered products of $Z_\alpha$ are self-adjoint and we assume that all Casimirs are also self-adjoint.

To obtain a full-dimensional (not upper-triangular) form of the matrix $\mathbb A$ let us consider adjoint action by any transport matrix:
\begin{theorem}\label{th:A-gen}\cite{ChSh2}
Any matrix $\mathbb A':=M_\gamma^{\text{T}} S^{\text{T}}\mathbb A S M_\gamma$ where $M_\gamma$ is a (transport) matrix satisfying commutation relations of Theorem \ref{th:MM} and commuting with $\mathbb A=M_1^{\text{T}} M_2$ satisfies the quantum reflection equation of Theorem \ref{th:A}.
\end{theorem}


\subsection{Casimirs of $b_n$-quiver}\label{ss:sln-quiver}
\begin{figure}[H]
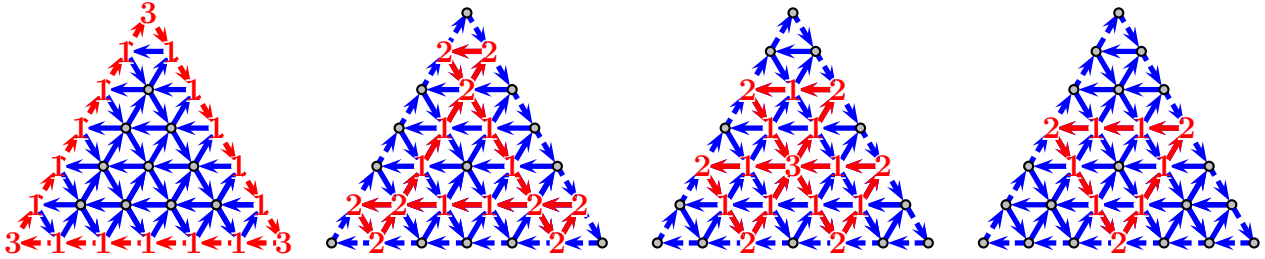


\caption{\small
Four central elements of the full-rank $b_6$-quiver.
}
\label{fi:Casimirs}
\end{figure}

For the full-rank $b_n$-quiver we have exactly $\bigl[\frac{n}{2}\bigr]+1$ Casimirs depicted in Fig~\ref{fi:Casimirs} for the example of $b_6$: numbers at vertices indicate the power with which the corresponding variables come into the product; all nonnumbered variables have power zero. All Casimirs correspond to closed broken-line paths in the $b_n$ quiver with reflections at the boundaries (the ``frozen'' variables at boundaries enter the product with powers two, powers of non-frozen variables can be 0,1,2, and 3, and they count how many times the path goes through the corresponding site). The total Poisson dimension of the full-rank quiver is therefore $\frac{(n+2)(n+1)}{2}-\bigl[\frac{n}{2}\bigr]-1$.

{\bf Notation.} Let 
\be\label{Ti}
T_i:=\col{\prod_{j=0}^{n-i} Z_{(j,n-j-i,i)}}
\ee
denote products of cluster variables along SE-diagonals of the $b_n$-quiver. 

\begin{figure}[H]
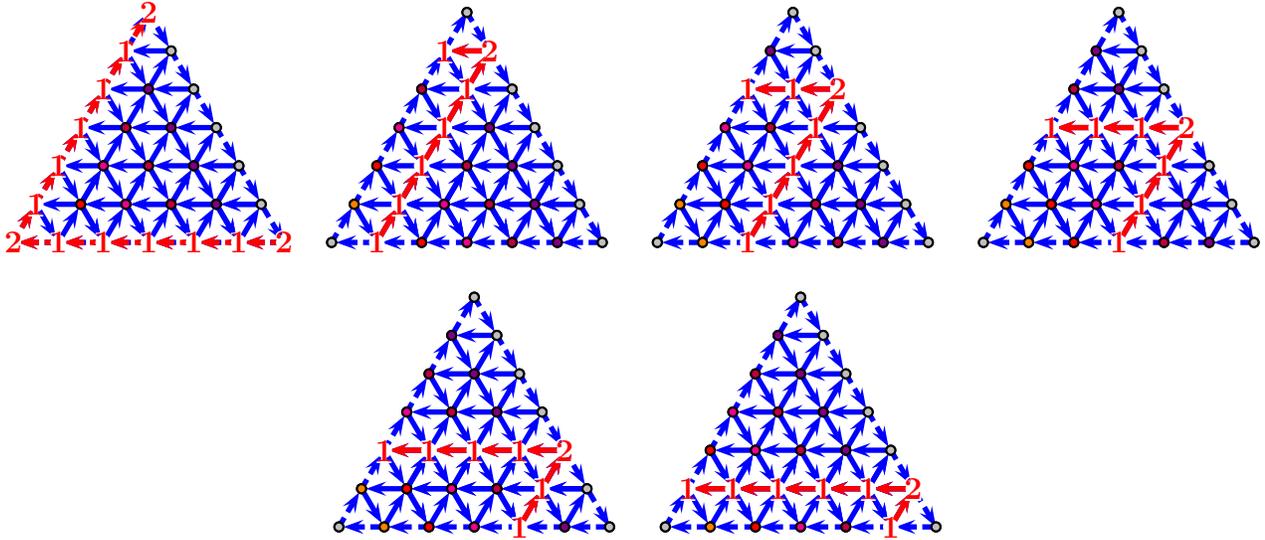

%
\caption{\small
Six mutually commuting elements $K_l$ (\ref{eq:Kl}) of the $b_6$-quiver; all these elements commute with all variables of the $\mathcal A_6$-quiver in Fig.~\ref{fi:amalgamation} and setting all $K_l=1$ results in the $\mathbb A^\hbar$ matrix of form (\ref{A-quantum}).
}
\label{fi:diagonal-Casimirs}
\end{figure}

\subsection{$\mathcal A_n$-quiver}\label{ss:An-quiver}

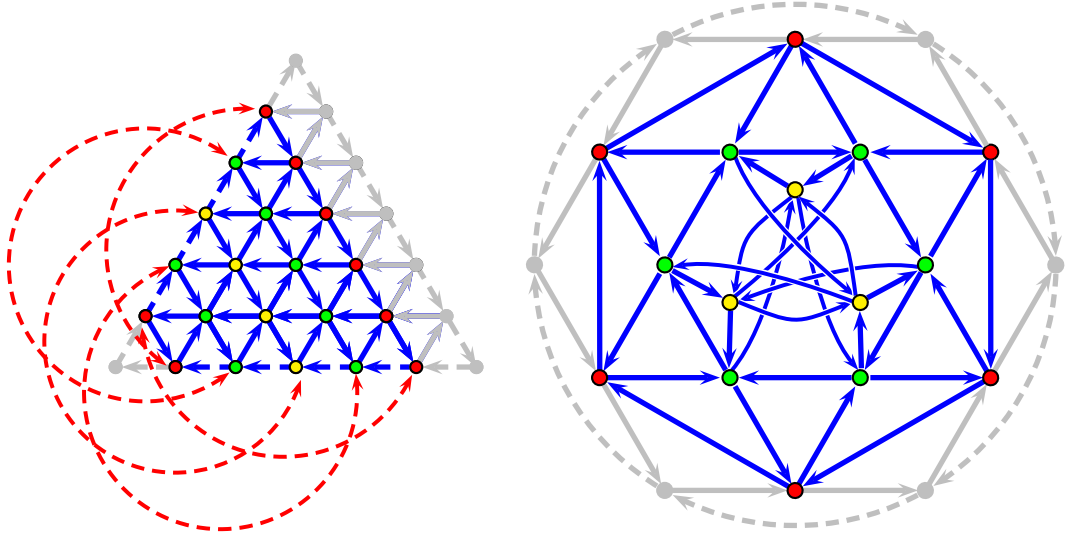
\begin{figure}[H]
\begin{pspicture}(-3,-3.5)(3,3.5){\psset{unit=0.8}
\newcommand{\PATGEN}{%
{\psset{unit=1}
\rput(0,0){\psline[linecolor=blue,linewidth=2pt]{->}(0,0)(.45,.765)}
\rput(0,0){\psline[linecolor=blue,linewidth=2pt]{->}(1,0)(0.1,0)}
\rput(0,0){\psline[linecolor=blue,linewidth=2pt]{->}(0,0)(.45,-.765)}
\put(0,0){\pscircle[fillstyle=solid,fillcolor=lightgray]{.1}}
}}
\newcommand{\PATLEFT}{%
{\psset{unit=1}
\rput(0,0){\psline[linecolor=blue,linewidth=2pt,linestyle=dashed]{->}(0,0)(.45,.765)}
\rput(0,0){\psline[linecolor=blue,linewidth=2pt]{->}(1,0)(0.1,0)}
\rput(0,0){\psline[linecolor=blue,linewidth=2pt]{->}(0,0)(.45,-.765)}
\put(0,0){\pscircle[fillstyle=solid,fillcolor=lightgray]{.1}}
}}
\newcommand{\PATRIGHT}{%
{\psset{unit=1}
\rput(0,0){\psline[linecolor=lightgray,linewidth=2pt,linestyle=dashed]{->}(0,0)(.45,-.765)}
\put(0,0){\pscircle[fillstyle=solid,fillcolor=lightgray]{.1}}
}}
\newcommand{\PATBOTTOM}{%
{\psset{unit=1}
\rput(0,0){\psline[linecolor=blue,linewidth=2pt]{->}(0,0)(.45,.765)}
\rput(0,0){\psline[linecolor=blue,linewidth=2pt,linestyle=dashed]{->}(1,0)(0.1,0)}
\put(0,0){\pscircle[fillstyle=solid,fillcolor=lightgray]{.1}}
}}
\newcommand{\PATTOP}{%
{\psset{unit=1}
\rput(0,0){\psline[linecolor=blue,linewidth=2pt]{->}(1,0)(0.1,0)}
\rput(0,0){\psline[linecolor=blue,linewidth=2pt]{->}(0,0)(.45,-.765)}
\put(0,0){\pscircle[fillstyle=solid,fillcolor=lightgray,linecolor=lightgray]{.1}}
}}
\newcommand{\PATBOTRIGHT}{%
{\psset{unit=1}
\rput(0,0){\psline[linecolor=blue,linewidth=2pt]{->}(0,0)(.45,.765)}
\put(0,0){\pscircle[fillstyle=solid,fillcolor=lightgray,linecolor=lightgray]{.1}}
\put(.5,0.85){\pscircle[fillstyle=solid,fillcolor=lightgray,linecolor=lightgray]{.1}}
}}
\newcommand{\PATNORTH}{%
			{\psset{unit=1}
				\rput(0,0){\psline[linecolor=lightgray,linewidth=2pt,linestyle=dashed]{->}(-.45,-.765)(0,0)}
				\rput(0,0){\psline[linecolor=lightgray,linewidth=2pt,linestyle=dashed]{<-}(.45,-.765)(0,0)}
				\put(0,0){\pscircle[fillstyle=solid,fillcolor=lightgray,linewidth=0pt]{.1}}
		}}
		\newcommand{\PATSW}{%
			{\psset{unit=1}
				\rput(0,0){\psline[linecolor=lightgray,linewidth=2pt,linestyle=dashed]{->}(0,0)(.45,.765)}
				\rput(0,0){\psline[linecolor=lightgray,linewidth=2pt,linestyle=dashed]{<-}(0.1,0)(0.9,0)}
				\put(0,0){\pscircle[fillstyle=solid,fillcolor=lightgray,linecolor=lightgray]{.1}}
		}}
		\newcommand{\PATSE}{%
			{\psset{unit=1}
				\rput(0,0){\psline[linecolor=lightgray,linewidth=2pt,linestyle=dashed]{<-}(0,0)(-.45,.765)}
				\rput(0,0){\psline[linecolor=lightgray,linewidth=2pt,linestyle=dashed]{<-}(-0.9,0)(-0.1,0)}
				\put(0,0){\pscircle[fillstyle=solid,fillcolor=lightgray,linecolor=lightgray]{.1}}
		}}
		
\multiput(-2.5,-0.85)(0.5,0.85){4}{\PATLEFT}
\multiput(-2,-1.7)(1,0){4}{\PATBOTTOM}
\put(-0.5,2.55){\PATTOP}
\multiput(-1.5,-0.85)(1,0){4}{\PATGEN}
\multiput(-1,0)(1,0){3}{\PATGEN}
\multiput(-.5,0.85)(1,0){2}{\PATGEN}
\put(0,1.7){\PATGEN}
\multiput(-1.5,-0.85)(1,0){4}{\PATGEN}
\multiput(0.5,2.55)(0.5,-0.85){4}{\PATRIGHT}
\put(2,-1.7){\PATBOTRIGHT}
\put(0,3.4){\PATNORTH}
\put(-3,-1.7){\PATSW}
\put(3,-1.7){\PATSE}
		\multiput(-0.5,2.55)(0.5,-0.85){5}{\rput(0,0){\psline[linecolor=lightgray,linewidth=2pt]{->}(1,0)(0.1,0)} }
		\multiput(0,1.7)(0.5,-0.85){5}{\rput(0,0){\psline[linecolor=lightgray,linewidth=2pt]{->}(0,0)(.45,.765)} }
		\multiput(0,3.4)(0.5,-0.85){7}{\pscircle[fillstyle=solid,fillcolor=lightgray,linecolor=lightgray]{.1}}
		\multiput(-0.5,2.55)(0.5,-0.85){6}{\pscircle[fillstyle=solid,fillcolor=red]{.1}}
		\multiput(-1,1.7)(0.5,-0.85){5}{\pscircle[fillstyle=solid,fillcolor=green]{.1}}
		\multiput(-1.5,0.85)(0.5,-0.85){4}{\pscircle[fillstyle=solid,fillcolor=yellow]{.1}}
		\multiput(-2,0)(0.5,-0.85){3}{\pscircle[fillstyle=solid,fillcolor=green]{.1}}
		\multiput(-2.5,-0.85)(0.5,-0.85){2}{\pscircle[fillstyle=solid,fillcolor=red]{.1}}
		\multiput(-3,-1.7)(0.5,-0.85){1}{\pscircle[fillstyle=solid,fillcolor=lightgray,linecolor=lightgray]{.1}}
\psarc[linecolor=red, linewidth=1.5pt,linestyle=dashed]{<->}(-2,-1.275){2.18}{80}{343}
\psarc[linecolor=red, linewidth=1.5pt,linestyle=dashed]{<->}(-2.5,0){2.27}{53}{308}
\psarc[linecolor=red, linewidth=1.5pt,linestyle=dashed]{<->}(-1.25,-2.125){2.27}{112}{368}
\psarc[linecolor=red, linewidth=1.5pt,linestyle=dashed]{<->}(-0.8,0.25){2.35}{87}{235}
\rput{100}(-0.1,0){\psarc[linecolor=red, linewidth=1.5pt,linestyle=dashed]{<->}(-0.8,0.25){2.36}{85}{235}}
%
}
\end{pspicture}
\begin{pspicture}(-3.5,-3.5)(3.5,3.5)
{\psset{unit=1}
\newcommand{\PATGEN}{%
{\psset{unit=1}
\psline[linecolor=blue,linewidth=2pt]{<-}(1,1.5)(2.5,1.5)
\psline[linecolor=blue,linewidth=2pt]{->}(-0.75,1.5)(0.75,1.5)
\psline[linecolor=blue,linewidth=2pt]{->}(-1,1.5)(-2.5,1.5)
\psline[linecolor=blue,linewidth=2pt]{->}(0.086,2.95)(2.5,1.56)
\pscircle[fillstyle=solid,fillcolor=green](0.866,1.5){.1}
\pscircle[fillstyle=solid,fillcolor=red](2.6,1.5){.1}
\psline[linecolor=lightgray,linewidth=2pt]{<-}(0.1,3)(1.634,3)
\psline[linecolor=lightgray,linewidth=2pt]{->}(-0.1,3)(-1.634,3)
\psarc[linecolor=lightgray, linewidth=2pt,linestyle=dashed]{<-}(0,0){3.468}{2}{58}
\pscircle[fillstyle=solid,fillcolor=lightgray,linecolor=lightgray](1.734,3){.1}
}}
\newcommand{\PATIN}{%
{\psset{unit=1}
\psline[linecolor=blue,linewidth=2pt]{<-}(-0.75,1.45)(-0.1,1.05)
\psline[linecolor=blue,linewidth=2pt]{->}(0.75,1.45)(0.1,1.05)
%
\psbezier[linecolor=white,linewidth=2.5pt]{->}(-0.84,-1.4)(-0.34,-0.8)(-0.1,0.5)(-0.04,0.9)
\psbezier[linecolor=blue,linewidth=1.5pt]{->}(-0.84,-1.4)(-0.34,-0.8)(-0.1,0.5)(-0.04,0.9)
\psbezier[linecolor=white,linewidth=2.5pt]{->}(-0.05,0.9)(-0.692,0.4)(-0.692,0.4)(-0.81,-0.42)
\psbezier[linecolor=blue,linewidth=1.5pt]{->}(-0.05,0.9)(-0.692,0.4)(-0.692,0.4)(-0.81,-0.42)
\psbezier[linecolor=white,linewidth=2.5pt]{<-}(0.84,-1.4)(0.34,-0.8)(0.1,0.5)(0.04,0.9)
\psbezier[linecolor=blue,linewidth=1.5pt]{<-}(0.84,-1.4)(0.34,-0.8)(0.1,0.5)(0.04,0.9)
\pscircle[fillstyle=solid,fillcolor=yellow](0,1){.1}
}}
\rput(0,0){\PATIN}
\rput{120}(0,0){\PATIN}
\rput{240}(0,0){\PATIN}
\rput(0,0){\PATGEN}
\rput{60}(0,0){\PATGEN}
\rput{120}(0,0){\PATGEN}
\rput{180}(0,0){\PATGEN}
\rput{240}(0,0){\PATGEN}
\rput{300}(0,0){\PATGEN}
}
\end{pspicture}
\caption{\small
The reduced quiver (the $\mathcal A_n$-quiver) obtained by amalgamating variables of the $b_n$-quiver and removing all remaining frozen variables.  (The example in the figure corresponds to $b_6$) Variables of the same color constitute Casimirs (with unit powers for all variables in every Casimir) depicted in Fig.~\ref{fi:Casimirs-A}.
}
\label{fi:amalgamation}
\end{figure}

We now construct the quiver corresponding to the system described by the reflection equation.  Note that transposition results, in particular, in amalgamation of variables $Z_{i,0,n-i}$ and $Z_{(0,n-i,i)}$ for $1\le i\le n-1$. We indicate these amalgamations by dashed arrows in Fig.~\ref{fi:amalgamation}. A more precise statement reads

\begin{lemma}\label{lm:amalgamation}
\begin{itemize}
\item matrix entries of $\mathbb A$ depend only on the variables $Z_{(i,j,k)}$ with $i,j,k>0$ and amalgamated variables $\col{Z_{(i,0,n-i)}Z_{(0,n-i,i)}}$, which are the variables of the {\em reduced quiver}, or $\mathcal A_n$-quiver, in Fig.~\ref{fi:amalgamation}, and on $n$ variables $K_l$, $1\le l\le n$, depicted in Fig.~\ref{fi:diagonal-Casimirs}:
\be\label{eq:Kl}
K_l=\col{\prod_{i=1}^l Z_{(i-1,n-l,l-i+1)} Z^2_{(l,n-l,0)}\prod_{j=1}^{n-l} Z_{(l,j-1,n-l-j+1)}  };
\ee 
\item elements $K_l$ commute mutually and commute with all variables of the reduced quiver;
\item upon setting all $K_l=1$ the matrix $\mathbb A$ takes the form (\ref{A-quantum}).
\end{itemize} 
\end{lemma}

We now unfreeze all amalgamated variables $\col{Z_{(i,0,n-i)}Z_{(0,n-i,i)}}$. From now on, we set all elements $K_l$ defined in (\ref{eq:Kl}) equal the unity. This eliminates the dependence on all remaining frozen variables and we remain with the reduced, or $\mathcal A_n$-quiver in Fig.~\ref{fi:amalgamation}. in which we indicate additional amalgamations by dashed lines.

\begin{figure}[H]
\begin{pspicture}(-2.7,-3)(2.7,2.5){\psset{unit=0.7}
\newcommand{\PATGEN}{%
{\psset{unit=1}
\rput(0,0){\psline[linecolor=blue,linewidth=2pt]{->}(0,0)(.45,.765)}
\rput(0,0){\psline[linecolor=blue,linewidth=2pt]{->}(1,0)(0.1,0)}
\rput(0,0){\psline[linecolor=blue,linewidth=2pt]{->}(0,0)(.45,-.765)}
\put(0,0){\pscircle[fillstyle=solid,fillcolor=lightgray]{.1}}
}}
\newcommand{\PATLEFT}{%
{\psset{unit=1}
\rput(0,0){\psline[linecolor=blue,linewidth=2pt,linestyle=dashed]{->}(0,0)(.45,.765)}
\rput(0,0){\psline[linecolor=blue,linewidth=2pt]{->}(1,0)(0.1,0)}
\rput(0,0){\psline[linecolor=blue,linewidth=2pt]{->}(0,0)(.45,-.765)}
\put(0,0){\pscircle[fillstyle=solid,fillcolor=lightgray]{.1}}
}}
\newcommand{\PATRIGHT}{%
{\psset{unit=1}
\rput(0,0){\psline[linecolor=lightgray,linewidth=2pt,linestyle=dashed]{->}(0,0)(.45,-.765)}
\put(0,0){\pscircle[fillstyle=solid,fillcolor=lightgray]{.1}}
}}
\newcommand{\PATBOTTOM}{%
{\psset{unit=1}
\rput(0,0){\psline[linecolor=blue,linewidth=2pt]{->}(0,0)(.45,.765)}
\rput(0,0){\psline[linecolor=blue,linewidth=2pt,linestyle=dashed]{->}(1,0)(0.1,0)}
\put(0,0){\pscircle[fillstyle=solid,fillcolor=lightgray]{.1}}
}}
\newcommand{\PATTOP}{%
{\psset{unit=1}
\rput(0,0){\psline[linecolor=blue,linewidth=2pt]{->}(1,0)(0.1,0)}
\rput(0,0){\psline[linecolor=blue,linewidth=2pt]{->}(0,0)(.45,-.765)}
\put(0,0){\pscircle[fillstyle=solid,fillcolor=lightgray]{.1}}
}}
\newcommand{\PATBOTRIGHT}{%
{\psset{unit=1}
\rput(0,0){\psline[linecolor=blue,linewidth=2pt]{->}(0,0)(.45,.765)}
\put(0,0){\pscircle[fillstyle=solid,fillcolor=lightgray]{.1}}
\put(.5,0.85){\pscircle[fillstyle=solid,fillcolor=lightgray]{.1}}
}}
\newcommand{\ODIN}{%
{\psset{unit=1}
\put(0,0){\pscircle[fillstyle=solid,fillcolor=white,linecolor=white]{.15}}
\put(0,0){\makebox(0,0)[cc]{\hbox{\tcw{\large$\mathbf 1$}}}}
\put(0,0){\makebox(0,0)[cc]{\hbox{\tcr{$\mathbf 1$}}}}
}}
\newcommand{\PATNORTH}{%
			{\psset{unit=1}
				\rput(0,0){\psline[linecolor=lightgray,linewidth=2pt,linestyle=dashed]{->}(-.45,-.765)(0,0)}
				\rput(0,0){\psline[linecolor=lightgray,linewidth=2pt,linestyle=dashed]{<-}(.45,-.765)(0,0)}
				\put(0,0){\pscircle[fillstyle=solid,fillcolor=lightgray]{.1}}
		}}
		\newcommand{\PATSW}{%
			{\psset{unit=1}
				\rput(0,0){\psline[linecolor=lightgray,linewidth=2pt,linestyle=dashed]{->}(0,0)(.45,.765)}
				\rput(0,0){\psline[linecolor=lightgray,linewidth=2pt,linestyle=dashed]{<-}(0.1,0)(0.9,0)}
				\put(0,0){\pscircle[fillstyle=solid,fillcolor=lightgray]{.1}}
		}}
		\newcommand{\PATSE}{%
			{\psset{unit=1}
				\rput(0,0){\psline[linecolor=lightgray,linewidth=2pt,linestyle=dashed]{<-}(0,0)(-.45,.765)}
				\rput(0,0){\psline[linecolor=lightgray,linewidth=2pt,linestyle=dashed]{<-}(-0.9,0)(-0.1,0)}
				\put(0,0){\pscircle[fillstyle=solid,fillcolor=lightgray]{.1}}
		}}
		
\multiput(-2.5,-0.85)(0.5,0.85){4}{\PATLEFT}
\multiput(-2,-1.7)(1,0){4}{\PATBOTTOM}
\put(-0.5,2.55){\PATTOP}
\multiput(-1.5,-0.85)(1,0){4}{\PATGEN}
\multiput(-1,0)(1,0){3}{\PATGEN}
\multiput(-.5,0.85)(1,0){2}{\PATGEN}
\put(0,1.7){\PATGEN}
\multiput(-1.5,-0.85)(1,0){4}{\PATGEN}
\multiput(0.5,2.55)(0.5,-0.85){4}{\PATRIGHT}
\put(2,-1.7){\PATBOTRIGHT}
\put(0,3.4){\PATNORTH}
\put(-3,-1.7){\PATSW}
\put(3,-1.7){\PATSE}
		\multiput(-0.5,2.55)(0.5,-0.85){5}{\rput(0,0){\psline[linecolor=lightgray,linewidth=2pt]{->}(1,0)(0.1,0)} }
		\multiput(0,1.7)(0.5,-0.85){5}{\rput(0,0){\psline[linecolor=lightgray,linewidth=2pt]{->}(0,0)(.45,.765)} }
		\multiput(0,3.4)(0.5,-0.85){7}{\pscircle[fillstyle=solid,fillcolor=lightgray,linecolor=lightgray]{.1}}
		\multiput(-0.5,2.55)(0.5,-0.85){6}{\pscircle[fillstyle=solid,fillcolor=violet]{.1}}
		\multiput(-1,1.7)(0.5,-0.85){5}{\pscircle[fillstyle=solid,fillcolor=purple]{.1}}
		\multiput(-1.5,0.85)(0.5,-0.85){4}{\pscircle[fillstyle=solid,fillcolor=magenta]{.1}}
		\multiput(-2,0)(0.5,-0.85){3}{\pscircle[fillstyle=solid,fillcolor=red]{.1}}
		\multiput(-2.5,-0.85)(0.5,-0.85){2}{\pscircle[fillstyle=solid,fillcolor=orange]{.1}}
		\multiput(-3,-1.7)(0.5,-0.85){1}{\pscircle[fillstyle=solid,fillcolor=lightgray,linecolor=lightgray]{.1}}
\psarc[linecolor=red, linewidth=1pt,linestyle=dashed]{<->}(-2,-1.275){2.18}{80}{343}
\psarc[linecolor=red, linewidth=1pt,linestyle=dashed]{<->}(-2.5,0){2.27}{53}{308}
\psarc[linecolor=red, linewidth=1pt,linestyle=dashed]{<->}(-1.25,-2.125){2.27}{112}{368}
\psarc[linecolor=red, linewidth=2pt,linestyle=dashed]{<->}(-0.8,0.25){2.35}{87}{235}
\rput{100}(-0.1,0){\psarc[linecolor=red, linewidth=2pt,linestyle=dashed]{<->}(-0.8,0.25){2.36}{85}{235}}
%
\multiput(-0.5,2.55)(0.5,-0.85){5}{\psline[linecolor=red,linewidth=2pt]{->}(0,0)(.45,-.765)}
\multiput(-2.5,-0.85)(0.5,-0.85){1}{\psline[linecolor=red,linewidth=2pt]{->}(0,0)(.45,-.765)}
\multiput(-.5,2.55)(0.5,-0.85){6}{\ODIN}
\multiput(-2.5,-0.85)(0.5,-0.85){2}{\ODIN}
%
}
\end{pspicture}
\begin{pspicture}(-2.7,-3)(2.7,2.5){\psset{unit=0.7}
\newcommand{\PATGEN}{%
{\psset{unit=1}
\rput(0,0){\psline[linecolor=blue,linewidth=2pt]{->}(0,0)(.45,.765)}
\rput(0,0){\psline[linecolor=blue,linewidth=2pt]{->}(1,0)(0.1,0)}
\rput(0,0){\psline[linecolor=blue,linewidth=2pt]{->}(0,0)(.45,-.765)}
\put(0,0){\pscircle[fillstyle=solid,fillcolor=lightgray]{.1}}
}}
\newcommand{\PATLEFT}{%
{\psset{unit=1}
\rput(0,0){\psline[linecolor=blue,linewidth=2pt,linestyle=dashed]{->}(0,0)(.45,.765)}
\rput(0,0){\psline[linecolor=blue,linewidth=2pt]{->}(1,0)(0.1,0)}
\rput(0,0){\psline[linecolor=blue,linewidth=2pt]{->}(0,0)(.45,-.765)}
\put(0,0){\pscircle[fillstyle=solid,fillcolor=lightgray]{.1}}
}}
\newcommand{\PATRIGHT}{%
{\psset{unit=1}
\rput(0,0){\psline[linecolor=lightgray,linewidth=2pt,linestyle=dashed]{->}(0,0)(.45,-.765)}
\put(0,0){\pscircle[fillstyle=solid,fillcolor=lightgray]{.1}}
}}
\newcommand{\PATBOTTOM}{%
{\psset{unit=1}
\rput(0,0){\psline[linecolor=blue,linewidth=2pt]{->}(0,0)(.45,.765)}
\rput(0,0){\psline[linecolor=blue,linewidth=2pt,linestyle=dashed]{->}(1,0)(0.1,0)}
\put(0,0){\pscircle[fillstyle=solid,fillcolor=lightgray]{.1}}
}}
\newcommand{\PATTOP}{%
{\psset{unit=1}
\rput(0,0){\psline[linecolor=blue,linewidth=2pt]{->}(1,0)(0.1,0)}
\rput(0,0){\psline[linecolor=blue,linewidth=2pt]{->}(0,0)(.45,-.765)}
\put(0,0){\pscircle[fillstyle=solid,fillcolor=lightgray]{.1}}
}}
\newcommand{\PATBOTRIGHT}{%
{\psset{unit=1}
\rput(0,0){\psline[linecolor=blue,linewidth=2pt]{->}(0,0)(.45,.765)}
\put(0,0){\pscircle[fillstyle=solid,fillcolor=lightgray]{.1}}
\put(.5,0.85){\pscircle[fillstyle=solid,fillcolor=lightgray]{.1}}
}}
\newcommand{\ODIN}{%
{\psset{unit=1}
\put(0,0){\pscircle[fillstyle=solid,fillcolor=white,linecolor=white]{.15}}
\put(0,0){\makebox(0,0)[cc]{\hbox{\tcw{\large$\mathbf 1$}}}}
\put(0,0){\makebox(0,0)[cc]{\hbox{\tcr{$\mathbf 1$}}}}
}}
\newcommand{\PATNORTH}{%
			{\psset{unit=1}
				\rput(0,0){\psline[linecolor=lightgray,linewidth=2pt,linestyle=dashed]{->}(-.45,-.765)(0,0)}
				\rput(0,0){\psline[linecolor=lightgray,linewidth=2pt,linestyle=dashed]{<-}(.45,-.765)(0,0)}
				\put(0,0){\pscircle[fillstyle=solid,fillcolor=lightgray]{.1}}
		}}
		\newcommand{\PATSW}{%
			{\psset{unit=1}
				\rput(0,0){\psline[linecolor=lightgray,linewidth=2pt,linestyle=dashed]{->}(0,0)(.45,.765)}
				\rput(0,0){\psline[linecolor=lightgray,linewidth=2pt,linestyle=dashed]{<-}(0.1,0)(0.9,0)}
				\put(0,0){\pscircle[fillstyle=solid,fillcolor=lightgray]{.1}}
		}}
		\newcommand{\PATSE}{%
			{\psset{unit=1}
				\rput(0,0){\psline[linecolor=lightgray,linewidth=2pt,linestyle=dashed]{<-}(0,0)(-.45,.765)}
				\rput(0,0){\psline[linecolor=lightgray,linewidth=2pt,linestyle=dashed]{<-}(-0.9,0)(-0.1,0)}
				\put(0,0){\pscircle[fillstyle=solid,fillcolor=lightgray]{.1}}
		}}
		
\multiput(-2.5,-0.85)(0.5,0.85){4}{\PATLEFT}
\multiput(-2,-1.7)(1,0){4}{\PATBOTTOM}
\put(-0.5,2.55){\PATTOP}
\multiput(-1.5,-0.85)(1,0){4}{\PATGEN}
\multiput(-1,0)(1,0){3}{\PATGEN}
\multiput(-.5,0.85)(1,0){2}{\PATGEN}
\put(0,1.7){\PATGEN}
\multiput(-1.5,-0.85)(1,0){4}{\PATGEN}
\multiput(0.5,2.55)(0.5,-0.85){4}{\PATRIGHT}
\put(2,-1.7){\PATBOTRIGHT}
\put(0,3.4){\PATNORTH}
\put(-3,-1.7){\PATSW}
\put(3,-1.7){\PATSE}
		\multiput(-0.5,2.55)(0.5,-0.85){5}{\rput(0,0){\psline[linecolor=lightgray,linewidth=2pt]{->}(1,0)(0.1,0)} }
		\multiput(0,1.7)(0.5,-0.85){5}{\rput(0,0){\psline[linecolor=lightgray,linewidth=2pt]{->}(0,0)(.45,.765)} }
		\multiput(0,3.4)(0.5,-0.85){7}{\pscircle[fillstyle=solid,fillcolor=lightgray,linecolor=lightgray]{.1}}
		\multiput(-0.5,2.55)(0.5,-0.85){6}{\pscircle[fillstyle=solid,fillcolor=violet]{.1}}
		\multiput(-1,1.7)(0.5,-0.85){5}{\pscircle[fillstyle=solid,fillcolor=purple]{.1}}
		\multiput(-1.5,0.85)(0.5,-0.85){4}{\pscircle[fillstyle=solid,fillcolor=magenta]{.1}}
		\multiput(-2,0)(0.5,-0.85){3}{\pscircle[fillstyle=solid,fillcolor=red]{.1}}
		\multiput(-2.5,-0.85)(0.5,-0.85){2}{\pscircle[fillstyle=solid,fillcolor=orange]{.1}}
		\multiput(-3,-1.7)(0.5,-0.85){1}{\pscircle[fillstyle=solid,fillcolor=lightgray,linecolor=lightgray]{.1}}
\psarc[linecolor=red, linewidth=1pt,linestyle=dashed]{<->}(-2,-1.275){2.18}{80}{343}
\psarc[linecolor=red, linewidth=2pt,linestyle=dashed]{<->}(-2.5,0){2.27}{53}{308}
\psarc[linecolor=red, linewidth=2pt,linestyle=dashed]{<->}(-1.25,-2.125){2.27}{112}{368}
\psarc[linecolor=red, linewidth=1pt,linestyle=dashed]{<->}(-0.8,0.25){2.35}{87}{235}
\rput{100}(-0.1,0){\psarc[linecolor=red, linewidth=1pt,linestyle=dashed]{<->}(-0.8,0.25){2.36}{85}{235}}
%
\multiput(-1,1.7)(0.5,-0.85){4}{\psline[linecolor=red,linewidth=2pt]{->}(0,0)(.45,-.765)}
\multiput(-2,0)(0.5,-0.85){2}{\psline[linecolor=red,linewidth=2pt]{->}(0,0)(.45,-.765)}
\multiput(-1,1.7)(0.5,-0.85){5}{\ODIN}
\multiput(-2,0)(0.5,-0.85){3}{\ODIN}
%
}
\end{pspicture}
\begin{pspicture}(-3,-3)(2,2.5){\psset{unit=0.7}
\newcommand{\PATGEN}{%
{\psset{unit=1}
\rput(0,0){\psline[linecolor=blue,linewidth=2pt]{->}(0,0)(.45,.765)}
\rput(0,0){\psline[linecolor=blue,linewidth=2pt]{->}(1,0)(0.1,0)}
\rput(0,0){\psline[linecolor=blue,linewidth=2pt]{->}(0,0)(.45,-.765)}
\put(0,0){\pscircle[fillstyle=solid,fillcolor=lightgray]{.1}}
}}
\newcommand{\PATLEFT}{%
{\psset{unit=1}
\rput(0,0){\psline[linecolor=blue,linewidth=2pt,linestyle=dashed]{->}(0,0)(.45,.765)}
\rput(0,0){\psline[linecolor=blue,linewidth=2pt]{->}(1,0)(0.1,0)}
\rput(0,0){\psline[linecolor=blue,linewidth=2pt]{->}(0,0)(.45,-.765)}
\put(0,0){\pscircle[fillstyle=solid,fillcolor=lightgray]{.1}}
}}
\newcommand{\PATRIGHT}{%
{\psset{unit=1}
\rput(0,0){\psline[linecolor=lightgray,linewidth=2pt,linestyle=dashed]{->}(0,0)(.45,-.765)}
\put(0,0){\pscircle[fillstyle=solid,fillcolor=lightgray]{.1}}
}}
\newcommand{\PATBOTTOM}{%
{\psset{unit=1}
\rput(0,0){\psline[linecolor=blue,linewidth=2pt]{->}(0,0)(.45,.765)}
\rput(0,0){\psline[linecolor=blue,linewidth=2pt,linestyle=dashed]{->}(1,0)(0.1,0)}
\put(0,0){\pscircle[fillstyle=solid,fillcolor=lightgray]{.1}}
}}
\newcommand{\PATTOP}{%
{\psset{unit=1}
\rput(0,0){\psline[linecolor=blue,linewidth=2pt]{->}(1,0)(0.1,0)}
\rput(0,0){\psline[linecolor=blue,linewidth=2pt]{->}(0,0)(.45,-.765)}
\put(0,0){\pscircle[fillstyle=solid,fillcolor=lightgray]{.1}}
}}
\newcommand{\PATBOTRIGHT}{%
{\psset{unit=1}
\rput(0,0){\psline[linecolor=blue,linewidth=2pt]{->}(0,0)(.45,.765)}
\put(0,0){\pscircle[fillstyle=solid,fillcolor=lightgray]{.1}}
\put(.5,0.85){\pscircle[fillstyle=solid,fillcolor=lightgray]{.1}}
}}
\newcommand{\ODIN}{%
{\psset{unit=1}
\put(0,0){\pscircle[fillstyle=solid,fillcolor=white,linecolor=white]{.15}}
\put(0,0){\makebox(0,0)[cc]{\hbox{\tcw{\large$\mathbf 1$}}}}
\put(0,0){\makebox(0,0)[cc]{\hbox{\tcr{$\mathbf 1$}}}}
}}
\newcommand{\PATNORTH}{%
			{\psset{unit=1}
				\rput(0,0){\psline[linecolor=lightgray,linewidth=2pt,linestyle=dashed]{->}(-.45,-.765)(0,0)}
				\rput(0,0){\psline[linecolor=lightgray,linewidth=2pt,linestyle=dashed]{<-}(.45,-.765)(0,0)}
				\put(0,0){\pscircle[fillstyle=solid,fillcolor=lightgray]{.1}}
		}}
		\newcommand{\PATSW}{%
			{\psset{unit=1}
				\rput(0,0){\psline[linecolor=lightgray,linewidth=2pt,linestyle=dashed]{->}(0,0)(.45,.765)}
				\rput(0,0){\psline[linecolor=lightgray,linewidth=2pt,linestyle=dashed]{<-}(0.1,0)(0.9,0)}
				\put(0,0){\pscircle[fillstyle=solid,fillcolor=lightgray]{.1}}
		}}
		\newcommand{\PATSE}{%
			{\psset{unit=1}
				\rput(0,0){\psline[linecolor=lightgray,linewidth=2pt,linestyle=dashed]{<-}(0,0)(-.45,.765)}
				\rput(0,0){\psline[linecolor=lightgray,linewidth=2pt,linestyle=dashed]{<-}(-0.9,0)(-0.1,0)}
				\put(0,0){\pscircle[fillstyle=solid,fillcolor=lightgray]{.1}}
		}}
		
\multiput(-2.5,-0.85)(0.5,0.85){4}{\PATLEFT}
\multiput(-2,-1.7)(1,0){4}{\PATBOTTOM}
\put(-0.5,2.55){\PATTOP}
\multiput(-1.5,-0.85)(1,0){4}{\PATGEN}
\multiput(-1,0)(1,0){3}{\PATGEN}
\multiput(-.5,0.85)(1,0){2}{\PATGEN}
\put(0,1.7){\PATGEN}
\multiput(-1.5,-0.85)(1,0){4}{\PATGEN}
\multiput(0.5,2.55)(0.5,-0.85){4}{\PATRIGHT}
\put(2,-1.7){\PATBOTRIGHT}
\put(0,3.4){\PATNORTH}
\put(-3,-1.7){\PATSW}
\put(3,-1.7){\PATSE}
		\multiput(-0.5,2.55)(0.5,-0.85){5}{\rput(0,0){\psline[linecolor=lightgray,linewidth=2pt]{->}(1,0)(0.1,0)} }
		\multiput(0,1.7)(0.5,-0.85){5}{\rput(0,0){\psline[linecolor=lightgray,linewidth=2pt]{->}(0,0)(.45,.765)} }
		\multiput(0,3.4)(0.5,-0.85){7}{\pscircle[fillstyle=solid,fillcolor=lightgray,linecolor=lightgray]{.1}}
		\multiput(-0.5,2.55)(0.5,-0.85){6}{\pscircle[fillstyle=solid,fillcolor=violet]{.1}}
		\multiput(-1,1.7)(0.5,-0.85){5}{\pscircle[fillstyle=solid,fillcolor=purple]{.1}}
		\multiput(-1.5,0.85)(0.5,-0.85){4}{\pscircle[fillstyle=solid,fillcolor=magenta]{.1}}
		\multiput(-2,0)(0.5,-0.85){3}{\pscircle[fillstyle=solid,fillcolor=red]{.1}}
		\multiput(-2.5,-0.85)(0.5,-0.85){2}{\pscircle[fillstyle=solid,fillcolor=orange]{.1}}
		\multiput(-3,-1.7)(0.5,-0.85){1}{\pscircle[fillstyle=solid,fillcolor=lightgray,linecolor=lightgray]{.1}}
\psarc[linecolor=red, linewidth=2pt,linestyle=dashed]{<->}(-2,-1.275){2.18}{80}{343}
\psarc[linecolor=red, linewidth=1pt,linestyle=dashed]{<->}(-2.5,0){2.27}{53}{308}
\psarc[linecolor=red, linewidth=1pt,linestyle=dashed]{<->}(-1.25,-2.125){2.27}{112}{368}
\psarc[linecolor=red, linewidth=1pt,linestyle=dashed]{<->}(-0.8,0.25){2.35}{87}{235}
\rput{100}(-0.1,0){\psarc[linecolor=red, linewidth=1pt,linestyle=dashed]{<->}(-0.8,0.25){2.36}{85}{235}}
%
\multiput(-1.5,0.85)(0.5,-0.85){3}{\psline[linecolor=red,linewidth=2pt]{->}(0,0)(.45,-.765)}
\multiput(-1.5,0.85)(0.5,-0.85){4}{\ODIN}
%
}
\end{pspicture}
\caption{\small
Three central elements for $\mathcal A_6$-quiver. 
}
\label{fi:Casimirs-A}
\end{figure}
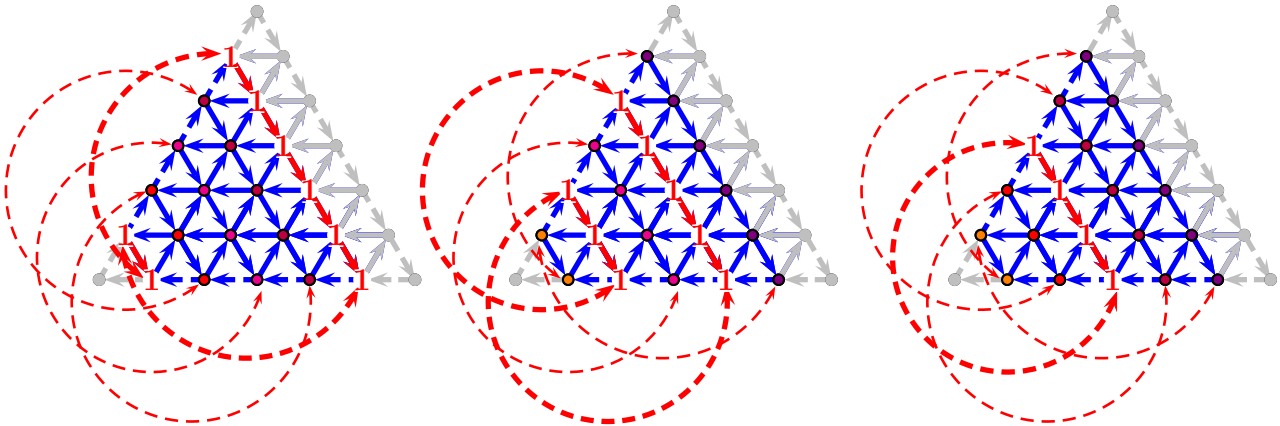

\begin{lemma}\label{lm:An-Casimirs}\cite{ChSh2}
Casimirs of $\mathcal A_n$-quiver are $\bigl[ \frac{n}{2}\bigr]$ elements $C_i$ given by the formula:
\be\label{Ci}
C_i=\col{T_i T_{n-i}}\hbox{ for }1\le i<n/2; \quad C_{n/2}=T_{n/2}\hbox{ if $n/2$ is an integer},
\ee
where $T_i$ are products (\ref{Ti}) of variables of $b_n$ quiver; note that $C_i$ depend only on amalgamated variables. The example for $n=6$ is depicted in Fig.~\ref{fi:Casimirs-A}. 
\end{lemma}

Since all Casimirs of $\mathcal A_n$ are generated by $\lambda$-power expansion terms for $\det(\mathbb A-\lambda \mathbb A^{\text{T}})$, we automatically obtain a semiclassical statement that $\det(\mathbb A-\lambda \mathbb A^{\text{T}})=P(C_1,\dots,C_{[n/2]})$, where $C_i$ are Casimirs of the $\mathcal A_n$-quiver. In the next section we explore this dependence.

\section{Solving characteristic equation}\label{s:character}
\setcounter{equation}{0}
The main theorem follows
\begin{theorem}\label{th:M}
Consider the $\mathcal A_n$-quiver in Fig.~\ref{fi:amalgamation}. Let $\lambda\in \mathbb C$ be an eigenvalue of the quantum operator $\mathbb A\mathbb A^{-\dagger}$, i.e., the number at which the equation $\bigl(\mathbb A\mathbb A^{-\dagger}-\lambda\mathbf I \bigr)\psi=0$ has a nontrivial null vector $\psi\in V\otimes W$, where $\mathbb A$ is $\mathbb A^\hbar$ from (\ref{A-quantum}). Then $n$ admissible values of $\lambda_i$ are
\be
\lambda_i=(-1)^{n-1} q^{-n}\times \left\{ \begin{array}{ll} 
\prod_{k=i}^{[n/2]} C_k&\hbox{for}\ 1\le i\le [n/2];\\ 1 &\hbox{for $i=(n+1)/2$ for odd $n$};\\ \prod_{k=n+1-i}^{[n/2]} C_k^{-1}&\hbox{for}\  n-[n/2]+1\le i \le n \end{array}\right.
\ee
\end{theorem}
{\bf Proof.} We use representation (\ref{eq:A-norm}) for the quantum matrix $\mathbb A^\hbar$. Since all entries of all matrices $\mathcal M_i$ are self-adjoint, $\mathcal M_i^\dagger=\mathcal M_i^{\text{T}}$, $Q^\dagger=Q^{-1}$ and $S^\dagger=S^{\text{T}}=(-1)^{n+1}S$. Then
\be\label{A-dag}
\bigl[\mathbb A^\hbar\bigr]^{\dagger}=\mathcal M_1^{\text{T}}S^{\text{T}}Q^{-1} \mathcal M_3^{\text{T}} \mathcal M_1,
\ee
so only the central block $S^{\text{T}}Q^{-1} \mathcal M_3^{\text{T}}$ sandwiched between $\mathcal M_1^{\text{T}}$ and $\mathcal M_1$ is changed, and
$$
\mathbb A^\hbar-\lambda \bigl[\mathbb A^\hbar\bigr]^{\dagger}=\mathcal M_1^{\text{T}} \bigl(\mathcal M_3QS  -\lambda S^{\text{T}}Q^{-1} \mathcal M_3^{\text{T}}\bigr)\mathcal M_1
$$
whereas the singularity equation becomes
\be\label{eq:sing}
 \bigl(\mathcal M_3QS  -\lambda S^{\text{T}}Q^{-1} \mathcal M_3^{\text{T}}\bigr)\psi=0
\ee
for some nonzero vector $\psi$. 

The {\em crucial observation} is that {\em both} matrices: $\mathcal M_3QS$ and  $S^{\text{T}}Q^{-1} \mathcal M_3^{\text{T}}$ in  (\ref{eq:sing}) {\em are upper-anti-diagonal matrices}! All solutions to the singularity equation therefore correspond to the situation where the combination of these matrices in (\ref{eq:sing}) has zero element on the anti-diagonal (then the determinant of this combination becomes zero). The matrix $\mathcal M_3$ is lower-triangular with the diagonal elements
$$
m_1=Z_{(n,0,0)},\quad m_i=\col {Z_{(n,0,0)}\prod_{j=1}^{i-1} T_i},\ 2\le i \le n,
$$
where $T_i$ are defined in (\ref{Ti}). The matrices 
$$
S^{\text{T}}Q^{-1} \mathcal M_3^{\text{T}}=\left[ \begin{array}{ccc} \star&\star& a_n\\ \star&\reflectbox{$\ddots$}&0\\a_1&0&0 \end{array}\right] \hbox{ and }  \mathcal M_3QS=\left[ \begin{array}{ccc} \star&\star& b_n\\ \star&\reflectbox{$\ddots$}&0\\b_1&0&0 \end{array}\right]
$$
then both are upper-anti-diagonal with
\be\label{ab-i}
a_i=(-1)^{i+1} q^{i-1/2} m_i,\quad b_i=(-1)^{n-i} q^{-n+i-1/2}m_{n+1-i},
\ee
and solutions of the singularity equation are
$$
\lambda_i=b_i/a_i=(-1)^{n-1}q^{-n} m_{n+1-i}(m_i)^{-1}
$$
Note that $m_i$ themselves are not Casimirs, but their quotients $m_{n+1-i}(m_i)^{-1}$ are: these quotients are just expressions in cases in Theorem~\ref{th:M}. This completes the proof. 

\begin{remark}\label{rm:general}
Note that representation (\ref{eq:A-norm}) remains valid for {\em any}, not necessarily lower-triangular matrix $\mathcal M_1$. This means that Theorem~\ref{th:M} remains valid in the case of general quantum $\mathbb A$-matrices, which may be not upper-triangular, provided they are presented in the cluster-variable decomposition form $B^\dagger \mathcal M_3QS B$ with $B$ having the same commutation relations with itself and with $\mathcal M_3$ as $\mathcal M_1$. In particular, it holds for any (amalgamated) $B=\mathcal M_1 QS \mathcal F$ with $\mathcal F$ commuting with all $\mathcal M_i$ and enjoying Lie--Poisson relations for its entries.
\end{remark}

The factors $(-1)^{n-1}q^{-n}$ in $\lambda_i$ in Theorem~\ref{th:M} are inessential and can be removed by a proper normalizaiton. A very interesting case is however when all $\lambda_i$ are different: in this case we have a full control over the Jordan form of $\mathbb A\mathbb A^{-\dagger}$, which becomes diagonal. Note also that since $\lambda_i$ are combinations of Casimirs, they are unit operators  $\mathbb I$ in the quantum space, so we conclude that 
$$
\mathbb A^\hbar = U^{-1} \Lambda U, \hbox{ where }\Lambda:=\sum_{i=1}^n \lambda_i e_{i,i}\otimes \mathbb I.
$$

\begin{corollary}\label{cor:thM}
Assuming all $\lambda_i$ are distinct and we have a natural $N$ such that  $\lambda_i^N=\lambda_j^N$ for all $i,j$ (so all $\lambda_i$ are proportional to distinct $N$th roots of the unity), we have that $\Bigl[ \mathbb A^\hbar \bigl[ \mathbb A^\hbar \bigr]^{-\dagger} \Bigr]^N=\hbox{const}\cdot \mathbf I$.
\end{corollary}

\section{$Sp_{2m}$ systems}\label{s:Sp2m}
\begin{figure}[h]
	\begin{pspicture}(-4,-2.5)(4,4){
		\newcommand{\LEFTDOWNARROW}{%
			{\psset{unit=1}
				\rput(0,0){\psline[doubleline=true,linewidth=1pt, doublesep=1pt, linecolor=black]{<-}(0,0)(.45,.23)}
		}}
		\newcommand{\DOWNARROW}{%
	{\psset{unit=1}
					\rput(0,0){\psline[doubleline=true,linewidth=1pt, doublesep=1pt, linecolor=black]{->}(0,0.1)(0,-0.566)}
		\put(0,0){\pscircle[fillstyle=solid,fillcolor=white]{.1}}
}}
		\newcommand{\LEFTUPARROW}{%
	{\psset{unit=1}
		\rput(0,0){\psline[doubleline=true,linewidth=1pt, doublesep=1pt, linecolor=black]{->}(0,0)(-.765,.45)}
}}
	\newcommand{\STARUP}{
			{\psset{unit=1}
	\rput(0,0){\psline[doubleline=true,linewidth=1pt, doublesep=1pt, linecolor=black]{<-}(0.06,-0.03)(.4,-.24)}
	\rput(0,0){\psline[doubleline=true,linewidth=1pt, doublesep=1pt, linecolor=black]{<-}(0,0.1)(0,.466)}
	\rput(0,0){\psline[doubleline=true,linewidth=1pt, doublesep=1pt, linecolor=black]{->}(0,0)(-.45,-.23)}
	\put(0,0){\pscircle[fillstyle=solid,fillcolor=black]{.1}}
	\put(0,.566){\pscircle[fillstyle=solid,fillcolor=white]{.1}}
}}
	\newcommand{\STARRIGNT}{
			{\psset{unit=1}
	\rput(0,0){\psline[doubleline=true,linewidth=1pt, doublesep=1pt, linecolor=black]{<-}(0.06,-0.03)(.4,-.24)}
	\rput(0,0){\psline[doubleline=true,linewidth=1pt, doublesep=1pt, linecolor=red]{->}(0,0.1)(0,.466)}
	\rput(0,0){\psline[doubleline=true,linewidth=1pt, doublesep=1pt, linecolor=lightgray]{->}(-0.5,0.95)(-0.5,1.25)}
	\rput(0,0){\psline[doubleline=true,linewidth=1pt, doublesep=1pt, linecolor=black]{<-}(0.07,0.6)(0.45,0.83)}
	\rput(0,0){\psline[doubleline=true,linewidth=1pt, doublesep=1pt, linecolor=black]{<-}(-0.43,0.8)(-0.08,0.616)}
	\put(0,0){\pscircle[fillstyle=solid,fillcolor=white]{.1}}
	\put(-0.5,0.85){\pscircle[fillstyle=solid,fillcolor=white]{.1}}
	\put(0,.566){\pscircle[fillstyle=solid,fillcolor=black]{.1}}
}}
		\newcommand{\PATGEN}{%
			{\psset{unit=1}
				\rput(0,0){\psline[linecolor=blue,linewidth=2pt]{->}(0,0)(.45,.765)}
				\rput(0,0){\psline[linecolor=blue,linewidth=2pt]{->}(1,0)(0.1,0)}
				\rput(0,0){\psline[linecolor=blue,linewidth=2pt]{->}(0,0)(.45,-.765)}
				\put(0,0){\pscircle[fillstyle=solid,fillcolor=red]{.1}}
		}}
		\newcommand{\PATLEFT}{%
			{\psset{unit=1}
				\rput(0,0){\psline[linecolor=blue,linewidth=2pt,linestyle=dashed]{->}(0,0)(.45,.765)}
				\rput(0,0){\psline[linecolor=blue,linewidth=2pt]{->}(1,0)(0.1,0)}
				\rput(0,0){\psline[linecolor=blue,linewidth=2pt]{->}(0,0)(.45,-.765)}
				\put(0,0){\pscircle[fillstyle=solid,fillcolor=red]{.1}}
		}}
		\newcommand{\PATRIGHT}{%
			{\psset{unit=1}
				\put(0,0){\pscircle[fillstyle=solid,fillcolor=lightgray]{.1}}
		}}
		\newcommand{\PATBOTTOM}{%
			{\psset{unit=1}
				\rput(0,0){\psline[linecolor=blue,linewidth=2pt]{->}(0,0)(.45,.765)}
				\rput(0,0){\psline[linecolor=blue,linewidth=2pt,linestyle=dashed]{->}(1,0)(0.1,0)}
				\put(0,0){\pscircle[fillstyle=solid,fillcolor=red]{.1}}
		}}
		\newcommand{\PATTOP}{%
			{\psset{unit=1}
				\rput(0,0){\psline[linecolor=blue,linewidth=2pt]{->}(1,0)(0.1,0)}
				\rput(0,0){\psline[linecolor=blue,linewidth=2pt]{->}(0,0)(.45,-.765)}
				\put(0,0){\pscircle[fillstyle=solid,fillcolor=red]{.1}}
		}}
		\newcommand{\PATBOTRIGHT}{%
			{\psset{unit=1}
				\rput(0,0){\psline[linecolor=blue,linewidth=2pt]{->}(0,0)(.45,.765)}
				\put(0,0){\pscircle[fillstyle=solid,fillcolor=red]{.1}}
				\put(.5,0.85){\pscircle[fillstyle=solid,fillcolor=lightgray]{.1}}
		}}
		\newcommand{\PATNORTH}{%
			{\psset{unit=1}
				\rput(0,0){\psline[linecolor=blue,linewidth=2pt,linestyle=dashed]{->}(-.45,-.765)(0,0)}
				\put(0,0){\pscircle[fillstyle=solid,fillcolor=lightgray]{.1}}
		}}
		\newcommand{\PATSW}{%
			{\psset{unit=1}
				\rput(0,0){\psline[linecolor=blue,linewidth=2pt,linestyle=dashed]{->}(0,0)(.45,.765)}
				\rput(0,0){\psline[linecolor=blue,linewidth=2pt,linestyle=dashed]{<-}(0.1,0)(0.9,0)}
				\put(0,0){\pscircle[fillstyle=solid,fillcolor=lightgray]{.1}}
		}}
		\newcommand{\PATSE}{%
			{\psset{unit=1}
				\rput(0,0){\psline[linecolor=blue,linewidth=2pt,linestyle=dashed]{<-}(-0.9,0)(-0.1,0)}
				\put(0,0){\pscircle[fillstyle=solid,fillcolor=lightgray]{.1}}
		}}
				\newcommand{\PATSPM}{%
			{\psset{unit=1}
				\rput(0,0){\psline[linecolor=blue,linewidth=2pt]{<-}(0.05,0.083)(.45,.765)}
				\rput(0,0){\psline[linecolor=blue,linewidth=2pt]{<-}(0.9,0)(0.1,0)}
				\rput(0,0){\psline[linecolor=blue,linestyle=dashed,linewidth=1.5pt]{->}(1,0)(.55,.765)}
				\rput(0,0){\psline[linecolor=blue,linestyle=dashed,linewidth=1.5pt]{->}(1,0)(.55,-.765)}
				\psarc[linecolor=black,linestyle=dashed,linewidth=1pt]{<->}(0,0){1}{65}{115}
				\put(0,0){\pscircle[fillstyle=solid,fillcolor=orange]{.1}}
				\put(0.5,0.85){\pscircle[fillstyle=solid,fillcolor=lightgray]{.1}}
				\put(1,0){\pscircle[fillstyle=solid,fillcolor=lightgray]{.1}}
				\put(1,0){\makebox(0,0)[cc]{\hbox{{\tiny$\times$}}}}
		}}
		\multiput(-2.5,-0.85)(0.5,0.85){4}{\PATLEFT}
		\multiput(-2,-1.7)(1,0){4}{\PATBOTTOM}
		\put(-0.5,2.55){\PATTOP}
		\put(0,3.4){\PATNORTH}
		\multiput(-1.5,-0.85)(1,0){4}{\PATGEN}
		\multiput(-1,0)(1,0){3}{\PATGEN}
		\multiput(-.5,0.85)(1,0){2}{\PATGEN}
		\put(0,1.7){\PATGEN}
		\multiput(-1.5,-0.85)(1,0){4}{\PATGEN}
		\multiput(0.5,2.55)(0.5,-0.85){4}{\PATRIGHT}
		\put(2,-1.7){\PATBOTRIGHT}
		\put(-3,-1.7){\PATSW}
		\put(3,-1.7){\PATSE}
		\multiput(0.5,2.55)(1,-1.7){3}{\PATSPM}
		\multiput(-2,-1.176)(1.0,0){5}{\STARUP}
		\multiput(-1.5,-0.335)(1.0,0){4}{\STARUP}
		\multiput(-1.0,0.5)(1.0,0){3}{\STARUP}
		\multiput(-.5,1.4)(1.0,0){2}{\STARUP}
		\put(0,2.3){\STARUP}
		\multiput(2.6,-1.4)(-0.5,.85){6}{\LEFTDOWNARROW}
		\multiput(-2.6,-1.4)(0.5,.85){6}{\LEFTUPARROW}
		\multiput(-2.5,-1.5)(1.0,0){6}{\DOWNARROW}
		\multiput(1,2.266)(1,-1.7){3}{\STARRIGNT}
		\multiput(-0.5,2.55)(0.5,-0.85){6}{\pscircle[fillstyle=solid,fillcolor=violet]{.1}}
		\multiput(-1,1.7)(0.5,-0.85){5}{\pscircle[fillstyle=solid,fillcolor=purple]{.1}}
		\multiput(-1.5,0.85)(0.5,-0.85){4}{\pscircle[fillstyle=solid,fillcolor=magenta]{.1}}
		\multiput(-2,0)(0.5,-0.85){3}{\pscircle[fillstyle=solid,fillcolor=red]{.1}}
		\multiput(-2.5,-0.85)(0.5,-0.85){2}{\pscircle[fillstyle=solid,fillcolor=orange]{.1}}
		\put(1.7,3.1){\makebox(0,0)[br]{\hbox{{$1$}}}}
		\put(1.7,2){\makebox(0,0)[br]{\hbox{{$2$}}}}
		\put(2.7,1.4){\makebox(0,0)[br]{\hbox{{$3$}}}}
		\put(2.7,0.3){\makebox(0,0)[br]{\hbox{{$4$}}}}
		\put(3.7,-0.3){\makebox(0,0)[br]{\hbox{{$5$}}}}
		\put(3.7,-1.4){\makebox(0,0)[br]{\hbox{{$6$}}}}

		\put(-0.9,3.2){\makebox(0,0)[br]{\hbox{{$1'$}}}}
		\put(-1.4,2.4){\makebox(0,0)[br]{\hbox{{$2'$}}}}
		\put(-1.9,1.6){\makebox(0,0)[br]{\hbox{{$3'$}}}}
		\put(-2.4,0.8){\makebox(0,0)[br]{\hbox{{$4'$}}}}
		\put(-2.9,-0.1){\makebox(0,0)[br]{\hbox{{$5'$}}}}
		\put(-3.4,-1.0){\makebox(0,0)[br]{\hbox{{$6'$}}}}
		
		\put(-2.2,-2.4){\makebox(0,0)[br]{\hbox{{$1''$}}}}
		\put(-1.2,-2.4){\makebox(0,0)[br]{\hbox{{$2''$}}}}
		\put(-0.2,-2.4){\makebox(0,0)[br]{\hbox{{$3''$}}}}
		\put(0.8,-2.4){\makebox(0,0)[br]{\hbox{{$4''$}}}}
		\put(1.8,-2.4){\makebox(0,0)[br]{\hbox{{$5''$}}}}
		\put(2.8,-2.4){\makebox(0,0)[br]{\hbox{{$6''$}}}}
	}
	\end{pspicture}
	\begin{pspicture}(-4,-2.5)(4,4){
		\newcommand{\LEFTDOWNARROW}{%
			{\psset{unit=1}
				\rput(0,0){\psline[doubleline=true,linewidth=1pt, doublesep=1pt, linecolor=black]{->}(0,0)(.765,.45)}
		}}
		\newcommand{\DOWNARROW}{%
	{\psset{unit=1}
					\rput(0,0){\psline[doubleline=true,linewidth=1pt, doublesep=1pt, linecolor=black]{->}(0,0.1)(0,-0.566)}
		\put(0,0){\pscircle[fillstyle=solid,fillcolor=white]{.1}}
}}
		\newcommand{\LEFTUPARROW}{%
	{\psset{unit=1}
		\rput(0,0){\psline[doubleline=true,linewidth=1pt, doublesep=1pt, linecolor=black]{<-}(0,0)(-.765,.45)}
}}
	\newcommand{\STARUP}{
			{\psset{unit=1}
	\rput(0,0){\psline[doubleline=true,linewidth=1pt, doublesep=1pt, linecolor=black]{->}(0.06,-0.03)(.4,-.24)}
	\rput(0,0){\psline[doubleline=true,linewidth=1pt, doublesep=1pt, linecolor=black]{<-}(0,0.1)(0,.466)}
	\rput(0,0){\psline[doubleline=true,linewidth=1pt, doublesep=1pt, linecolor=black]{<-}(0,0)(-.45,-.23)}
	\put(0,0){\pscircle[fillstyle=solid,fillcolor=black]{.1}}
	\put(0,.566){\pscircle[fillstyle=solid,fillcolor=white]{.1}}
}}
		\newcommand{\PATGEN}{%
			{\psset{unit=1}
				\rput(0,0){\psline[linecolor=blue,linewidth=2pt]{->}(0,0)(.45,.765)}
				\rput(0,0){\psline[linecolor=blue,linewidth=2pt]{->}(1,0)(0.1,0)}
				\rput(0,0){\psline[linecolor=blue,linewidth=2pt]{->}(0,0)(.45,-.765)}
				\put(0,0){\pscircle[fillstyle=solid,fillcolor=red]{.1}}
		}}
		\newcommand{\PATLEFT}{%
			{\psset{unit=1}
				\rput(0,0){\psline[linecolor=blue,linewidth=2pt,linestyle=dashed]{->}(0,0)(.45,.765)}
				\rput(0,0){\psline[linecolor=blue,linewidth=2pt]{->}(1,0)(0.1,0)}
				\rput(0,0){\psline[linecolor=blue,linewidth=2pt]{->}(0,0)(.45,-.765)}
				\put(0,0){\pscircle[fillstyle=solid,fillcolor=red]{.1}}
		}}
		\newcommand{\PATRIGHT}{%
			{\psset{unit=1}
				\put(0,0){\pscircle[fillstyle=solid,fillcolor=lightgray]{.1}}
		}}
		\newcommand{\PATBOTTOM}{%
			{\psset{unit=1}
				\rput(0,0){\psline[linecolor=blue,linewidth=2pt]{->}(0,0)(.45,.765)}
				\rput(0,0){\psline[linecolor=blue,linewidth=2pt,linestyle=dashed]{->}(1,0)(0.1,0)}
				\put(0,0){\pscircle[fillstyle=solid,fillcolor=red]{.1}}
		}}
		\newcommand{\PATTOP}{%
			{\psset{unit=1}
				\rput(0,0){\psline[linecolor=blue,linewidth=2pt]{->}(1,0)(0.1,0)}
				\rput(0,0){\psline[linecolor=blue,linewidth=2pt]{->}(0,0)(.45,-.765)}
				\put(0,0){\pscircle[fillstyle=solid,fillcolor=red]{.1}}
		}}
		\newcommand{\PATBOTRIGHT}{%
			{\psset{unit=1}
				\rput(0,0){\psline[linecolor=blue,linewidth=2pt]{->}(0,0)(.45,.765)}
				\put(0,0){\pscircle[fillstyle=solid,fillcolor=red]{.1}}
				\put(.5,0.85){\pscircle[fillstyle=solid,fillcolor=lightgray]{.1}}
		}}
		\newcommand{\PATNORTH}{%
			{\psset{unit=1}
				\rput(0,0){\psline[linecolor=blue,linewidth=2pt,linestyle=dashed]{->}(-.45,-.765)(0,0)}
				\put(0,0){\pscircle[fillstyle=solid,fillcolor=lightgray]{.1}}
		}}
		\newcommand{\PATSW}{%
			{\psset{unit=1}
				\rput(0,0){\psline[linecolor=blue,linewidth=2pt,linestyle=dashed]{->}(0,0)(.45,.765)}
				\rput(0,0){\psline[linecolor=blue,linewidth=2pt,linestyle=dashed]{<-}(0.1,0)(0.9,0)}
				\put(0,0){\pscircle[fillstyle=solid,fillcolor=lightgray]{.1}}
		}}
		\newcommand{\PATSE}{%
			{\psset{unit=1}
				\rput(0,0){\psline[linecolor=blue,linewidth=2pt,linestyle=dashed]{<-}(-0.9,0)(-0.1,0)}
				\put(0,0){\pscircle[fillstyle=solid,fillcolor=lightgray]{.1}}
		}}
					\newcommand{\PATS}{%
			{\psset{unit=1}
				\rput(0,0){\psline[linecolor=blue,linewidth=2pt]{<-}(-0.05,0.083)(-.45,.765)}
				\rput(0,0){\psline[linecolor=blue,linewidth=2pt]{<-}(0.9,0)(0.1,0)}
				\rput(0,0){\psline[linecolor=blue,linestyle=dashed,linewidth=2pt]{->}(1,0)(-.4,.83)}
				\rput(0,0){\psline[linecolor=blue,linestyle=dashed,linewidth=2pt]{->}(1,0)(.55,-.765)}
				\put(0,0){\pscircle[fillstyle=solid,fillcolor=orange]{.1}}
				\put(1,0){\pscircle[fillstyle=solid,fillcolor=lightgray]{.1}}
				\put(1,0){\makebox(0,0)[cc]{\hbox{{\tiny$\times$}}}}
		}}
		\multiput(-2.5,-0.85)(0.5,0.85){4}{\PATLEFT}
		\multiput(-2,-1.7)(1,0){4}{\PATBOTTOM}
		\put(-0.5,2.55){\PATTOP}
		\put(0,3.4){\PATNORTH}
		\multiput(-1.5,-0.85)(1,0){4}{\PATGEN}
		\multiput(-1,0)(1,0){3}{\PATGEN}
		\multiput(-.5,0.85)(1,0){2}{\PATGEN}
		\put(0,1.7){\PATGEN}
		\multiput(-1.5,-0.85)(1,0){4}{\PATGEN}
		\multiput(0.5,2.55)(0.5,-0.85){4}{\PATRIGHT}
		\put(2,-1.7){\PATBOTRIGHT}
		\put(-3,-1.7){\PATSW}
		\put(3,-1.7){\PATSE}
		\multiput(0.5,2.55)(1,-1.7){3}{\PATS}
		%
		\multiput(-0.5,2.55)(0.5,-0.85){6}{\pscircle[fillstyle=solid,fillcolor=red]{.1}}
		\multiput(-1,1.7)(0.5,-0.85){5}{\pscircle[fillstyle=solid,fillcolor=green]{.1}}
		\multiput(-1.5,0.85)(0.5,-0.85){4}{\pscircle[fillstyle=solid,fillcolor=yellow]{.1}}
		\multiput(-2,0)(0.5,-0.85){3}{\pscircle[fillstyle=solid,fillcolor=green]{.1}}
		\multiput(-2.5,-0.85)(0.5,-0.85){2}{\pscircle[fillstyle=solid,fillcolor=red]{.1}}
		}
	\end{pspicture}
\caption{\small
	In the left picture we present the $Sp_{2m}$ directed network for $m=3$. We remove three sinks indicated in the light color and located at the NE side of the triangle and amalgamate cluster variables of faces separated by the corresponding directed edges (these amalgamations are denoted by dashed arcs). Note the appearance of  (three) double arrows directed upward: the transport matrix ceases to be upper triangular and becomes only block-upper triangular. The resulting $Sp_{2m}$-quiver is presented in the right picture.  
	}
	\label{fi:Sp2m}
\end{figure}
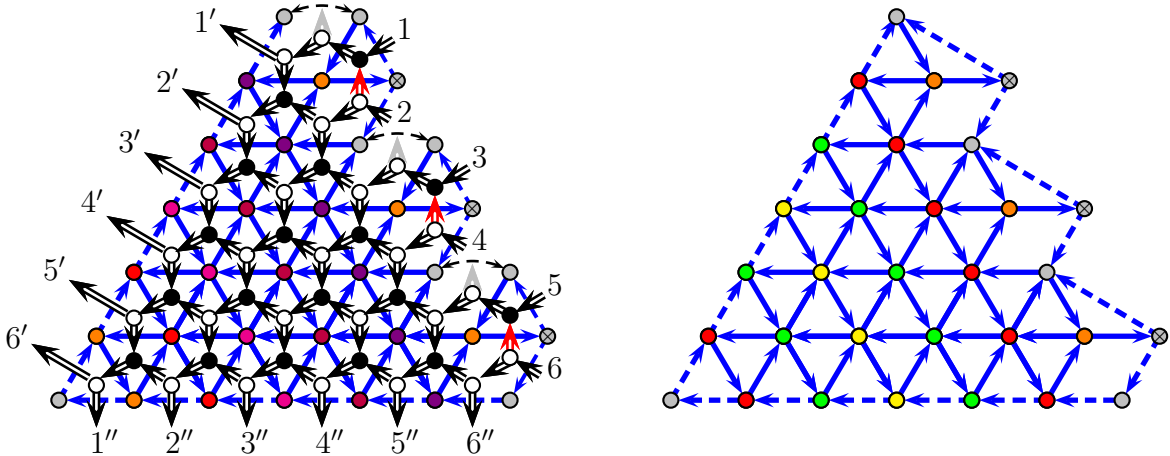

\setcounter{equation}{0}
In this section, we find Casimirs for the matrix $\mathbb A$ having $2\times 2$-block matrix form. These systems correspond to $Sp_{2m}$ twisted Yangians (see \cite{MRS}). Transport matrices for this system must have block-triangular form. A directed network corresponding to the case of $2\times 2$ blocks is presented in the left part of Fig.~\ref{fi:Sp2m}. Note that we have to eliminate $m$ redundant sinks in this picture and amalgamate variables of faces separated by the corresponding sinks. The resulting quiver is presented in the right part of the figure. 

\begin{figure}[h]
	\begin{pspicture}(-3.5,-2.5)(3.5,4)
		\newcommand{\PATGEN}{%
			{\psset{unit=1}
				\rput(0,0){\psline[linecolor=blue,linewidth=2pt]{->}(0,0)(.45,.765)}
				\rput(0,0){\psline[linecolor=blue,linewidth=2pt]{->}(1,0)(0.1,0)}
				\rput(0,0){\psline[linecolor=blue,linewidth=2pt]{->}(0,0)(.45,-.765)}
				\put(0,0){\pscircle[fillstyle=solid,fillcolor=red]{.1}}
		}}
		\newcommand{\PATLEFT}{%
			{\psset{unit=1}
				\rput(0,0){\psline[linecolor=blue,linewidth=2pt,linestyle=dashed]{->}(0,0)(.45,.765)}
				\rput(0,0){\psline[linecolor=blue,linewidth=2pt]{->}(1,0)(0.1,0)}
				\rput(0,0){\psline[linecolor=blue,linewidth=2pt]{->}(0,0)(.45,-.765)}
				\put(0,0){\pscircle[fillstyle=solid,fillcolor=red]{.1}}
		}}
		\newcommand{\PATRIGHT}{%
			{\psset{unit=1}
				\put(0,0){\pscircle[fillstyle=solid,fillcolor=lightgray]{.1}}
		}}
		\newcommand{\PATBOTTOM}{%
			{\psset{unit=1}
				\rput(0,0){\psline[linecolor=blue,linewidth=2pt]{->}(0,0)(.45,.765)}
				\rput(0,0){\psline[linecolor=blue,linewidth=2pt,linestyle=dashed]{->}(1,0)(0.1,0)}
				\put(0,0){\pscircle[fillstyle=solid,fillcolor=red]{.1}}
		}}
		\newcommand{\PATTOP}{%
			{\psset{unit=1}
				\rput(0,0){\psline[linecolor=blue,linewidth=2pt]{->}(1,0)(0.1,0)}
				\rput(0,0){\psline[linecolor=blue,linewidth=2pt]{->}(0,0)(.45,-.765)}
				\put(0,0){\pscircle[fillstyle=solid,fillcolor=red]{.1}}
		}}
		\newcommand{\PATBOTRIGHT}{%
			{\psset{unit=1}
				\rput(0,0){\psline[linecolor=blue,linewidth=2pt]{->}(0,0)(.45,.765)}
				\put(0,0){\pscircle[fillstyle=solid,fillcolor=red]{.1}}
				\put(.5,0.85){\pscircle[fillstyle=solid,fillcolor=lightgray]{.1}}
		}}
		\newcommand{\PATNORTH}{%
			{\psset{unit=1}
				\rput(0,0){\psline[linecolor=blue,linewidth=2pt,linestyle=dashed]{->}(-.45,-.765)(0,0)}
				\put(0,0){\pscircle[fillstyle=solid,fillcolor=lightgray]{.1}}
		}}
		\newcommand{\PATSW}{%
			{\psset{unit=1}
				\rput(0,0){\psline[linecolor=blue,linewidth=2pt,linestyle=dashed]{->}(0,0)(.45,.765)}
				\rput(0,0){\psline[linecolor=blue,linewidth=2pt,linestyle=dashed]{<-}(0.1,0)(0.9,0)}
				\put(0,0){\pscircle[fillstyle=solid,fillcolor=lightgray]{.1}}
		}}
		\newcommand{\PATSE}{%
			{\psset{unit=1}
				\rput(0,0){\psline[linecolor=blue,linewidth=2pt,linestyle=dashed]{<-}(-0.9,0)(-0.1,0)}
				\put(0,0){\pscircle[fillstyle=solid,fillcolor=lightgray]{.1}}
		}}
					\newcommand{\PATS}{%
			{\psset{unit=1}
				\rput(0,0){\psline[linecolor=blue,linewidth=2pt]{<-}(-0.05,0.083)(-.45,.765)}
				\rput(0,0){\psline[linecolor=blue,linewidth=2pt]{<-}(0.9,0)(0.1,0)}
				\rput(0,0){\psline[linecolor=blue,linestyle=dashed,linewidth=2pt]{->}(1,0)(-.4,.83)}
				\rput(0,0){\psline[linecolor=blue,linestyle=dashed,linewidth=2pt]{->}(1,0)(.55,-.765)}
				\put(0,0){\pscircle[fillstyle=solid,fillcolor=orange]{.1}}
				\put(1,0){\pscircle[fillstyle=solid,fillcolor=lightgray]{.1}}
				\put(1,0){\makebox(0,0)[cc]{\hbox{{\tiny$\times$}}}}
		}}
		\multiput(-2.5,-0.85)(0.5,0.85){4}{\PATLEFT}
		\multiput(-2,-1.7)(1,0){4}{\PATBOTTOM}
		\put(-0.5,2.55){\PATTOP}
		\put(0,3.4){\PATNORTH}
		\multiput(-1.5,-0.85)(1,0){4}{\PATGEN}
		\multiput(-1,0)(1,0){3}{\PATGEN}
		\multiput(-.5,0.85)(1,0){2}{\PATGEN}
		\put(0,1.7){\PATGEN}
		\multiput(-1.5,-0.85)(1,0){4}{\PATGEN}
		\multiput(0.5,2.55)(0.5,-0.85){4}{\PATRIGHT}
		\put(2,-1.7){\PATBOTRIGHT}
		\put(-3,-1.7){\PATSW}
		\put(3,-1.7){\PATSE}
		\multiput(0.5,2.55)(1,-1.7){3}{\PATS}
		\multiput(0,-1.7)(0.5,0.85){3}{\psline[linecolor=blue,linewidth=4pt]{->}(0,0)(.45,.765)\psline[linecolor=white,linewidth=2pt]{->}(0,0)(.4,.68)}
		\multiput(-1.5,0.85)(1,0){3}{\psline[linecolor=blue,linewidth=4pt]{<-}(0.1,0)(0.9,0)\psline[linecolor=white,linewidth=2pt]{<-}(0.2,0)(0.9,0)}

		\multiput(-0.5,2.55)(0.5,-0.85){6}{\pscircle[fillstyle=solid,fillcolor=red]{.1}}
		\multiput(-1,1.7)(0.5,-0.85){5}{\pscircle[fillstyle=solid,fillcolor=green]{.1}}
		\multiput(-1.5,0.85)(0.5,-0.85){4}{\pscircle[fillstyle=solid,fillcolor=yellow]{.1}}
		\multiput(-2,0)(0.5,-0.85){3}{\pscircle[fillstyle=solid,fillcolor=green]{.1}}
		\multiput(-2.5,-0.85)(0.5,-0.85){2}{\pscircle[fillstyle=solid,fillcolor=red]{.1}}
		\put(-3.1,-1.9){\makebox(0,0)[tr]{\hbox{{$S_1$}}}}
		\put(3.1,-1.9){\makebox(0,0)[tl]{\hbox{{$S_1$}}}}
		\put(0,3.6){\makebox(0,0)[bc]{\hbox{{$S_1$}}}}
		\put(0.5,2.4){\makebox(0,0)[tl]{\hbox{{$S_2$}}}}
		\put(1.2,1.65){\makebox(0,0)[bl]{\hbox{{$S_3$}}}}
		\put(1.5,0.7){\makebox(0,0)[tl]{\hbox{{$S_4$}}}}
		\put(2.2,-0.05){\makebox(0,0)[bl]{\hbox{{$S_5$}}}}
		\put(2.5,-1){\makebox(0,0)[tl]{\hbox{{$S_6$}}}}
		\put(1.5,2.75){\makebox(0,0)[bc]{\hbox{{$Q_2$}}}}
		\put(2.5,1.05){\makebox(0,0)[bc]{\hbox{{$Q_4$}}}}
		\put(3.5,-.65){\makebox(0,0)[bc]{\hbox{{$Q_6$}}}}
		\put(-0.6,2.65){\makebox(0,0)[br]{\hbox{{$a$}}}}
		\put(-1.1,1.8){\makebox(0,0)[br]{\hbox{{$b$}}}}
		\put(-1.6,0.95){\makebox(0,0)[br]{\hbox{{$c$}}}}
		\put(-2.1,0.1){\makebox(0,0)[br]{\hbox{{$d$}}}}
		\put(-2.6,-0.75){\makebox(0,0)[br]{\hbox{{$e$}}}}
		\put(-2,-2.2){\makebox(0,0)[bc]{\hbox{{$a$}}}}
		\put(-1,-2.2){\makebox(0,0)[bc]{\hbox{{$b$}}}}
		\put(0,-2.2){\makebox(0,0)[bc]{\hbox{{$c$}}}}
		\put(1,-2.2){\makebox(0,0)[bc]{\hbox{{$d$}}}}
		\put(2,-2.2){\makebox(0,0)[bc]{\hbox{{$e$}}}}
		\put(0.5,-1.45){\makebox(0,0)[bc]{\hbox{{$f$}}}}
		\put(1,-0.6){\makebox(0,0)[bc]{\hbox{{$g$}}}}
		\put(-0.47,1.15){\makebox(0,0)[bc]{\hbox{{$k$}}}}
		\put(0.53,1.15){\makebox(0,0)[bc]{\hbox{{$h$}}}}
	\end{pspicture}
\begin{pspicture}(-4,-3.3)(4,3.3)
{\psset{unit=1}
\newcommand{\PATGEN}{%
{\psset{unit=1}
\psline[linecolor=blue,linewidth=2pt]{<-}(1,1.5)(2.5,1.5)
\psline[linecolor=blue,linewidth=2pt]{->}(-0.75,1.5)(0.75,1.5)
\psline[linecolor=blue,linewidth=2pt]{->}(-1,1.5)(-2.5,1.5)
\psline[linecolor=blue,linewidth=2pt]{->}(0.1,2.94)(2.5,1.56)
\pscircle[fillstyle=solid,fillcolor=green](0.866,1.5){.1}
\pscircle[fillstyle=solid,fillcolor=red](2.6,1.5){.1}
}}
\newcommand{\PATIN}{%
{\psset{unit=1}
\psline[linecolor=blue,linewidth=2pt]{<-}(-0.75,1.45)(-0.1,1.05)
\psline[linecolor=blue,linewidth=2pt]{->}(0.75,1.45)(0.1,1.05)
%
\psbezier[linecolor=white,linewidth=2.5pt]{->}(-0.84,-1.4)(-0.34,-0.8)(-0.1,0.5)(-0.04,0.9)
\psbezier[linecolor=blue,linewidth=1.5pt]{->}(-0.84,-1.4)(-0.34,-0.8)(-0.1,0.5)(-0.04,0.9)
\psbezier[linecolor=white,linewidth=2.5pt]{->}(-0.05,0.9)(-0.692,0.4)(-0.692,0.4)(-0.81,-0.42)
\psbezier[linecolor=blue,linewidth=1.5pt]{->}(-0.05,0.9)(-0.692,0.4)(-0.692,0.4)(-0.81,-0.42)
\psbezier[linecolor=white,linewidth=2.5pt]{<-}(0.84,-1.4)(0.34,-0.8)(0.1,0.5)(0.04,0.9)
\psbezier[linecolor=blue,linewidth=1.5pt]{<-}(0.84,-1.4)(0.34,-0.8)(0.1,0.5)(0.04,0.9)
\pscircle[fillstyle=solid,fillcolor=yellow](0,1){.1}
}}
\newcommand{\PATOUT}{%
{\psset{unit=1}
\psline[linecolor=blue,linewidth=2pt]{<-}(2.65,1.4)(3.17,0.1)
\psline[linecolor=blue,linewidth=2pt]{->}(2.65,-1.4)(3.17,-0.1)
\pscircle[fillstyle=solid,fillcolor=orange](3.2,0){.1}
\psline[linecolor=blue,linewidth=2pt]{->}(3.3,0)(3.7,0)
\pscircle[fillstyle=solid,fillcolor=lightgray](3.8,0){.1}
\put(3.8,0){\makebox(0,0)[cc]{\hbox{{\tiny$\times$}}}}
%
\psbezier[linecolor=blue,linewidth=1.5pt,linestyle=dashed]{->}(3.8,0.1)(3.8,1.2)(2.7,2.77)(1.7,2.77)
\psbezier[linecolor=blue,linewidth=1.5pt,linestyle=dashed]{->}(3.8,-0.1)(3.8,-1.2)(2.7,-2.77)(1.7,-2.77)
\psbezier[linecolor=blue,linewidth=1.5pt]{<-}(3.25,0.1)(3.2,1.2)(2.7,2.33)(1.7,2.73)
}}
\newcommand{\PATOUTONE}{%
{\psset{unit=1}
\psline[linecolor=blue,linewidth=2pt]{<-}(2.65,1.4)(3.17,0.1)
\psline[linecolor=blue,linewidth=2pt]{->}(2.65,-1.4)(3.17,-0.1)
\pscircle[fillstyle=solid,fillcolor=lightgray](3.2,0){.1}
}}
\rput(0,0){\PATIN}
\rput{120}(0,0){\PATIN}
\rput{240}(0,0){\PATIN}
\rput(0,0){\PATGEN}
\rput{60}(0,0){\PATGEN}
\rput{120}(0,0){\PATGEN}
\rput{180}(0,0){\PATGEN}
\rput{240}(0,0){\PATGEN}
\rput{300}(0,0){\PATGEN}
\rput(0,0){\PATOUT}
\rput{60}(0,0){\PATOUTONE}
\rput{120}(0,0){\PATOUT}
\rput{180}(0,0){\PATOUTONE}
\rput{240}(0,0){\PATOUT}
\rput{300}(0,0){\PATOUTONE}
}
\end{pspicture}
\caption{\small
The $\mathcal A_{Sp_{2m}}$-quiver corresponding to block-upper triangular $\mathbb A$ in the case of $2\times 2$ blocks. The example in the figure corresponds to $m=3$. Amalgamated variables are indicated by the same letters $a$, $b$, $c$, $d$, $e$, and $S_1$ in the left figure. In the same figure we indicate cluster variables $c$, $f$, $g$, $S_4$, $h$, $k$ stretched along the path corresponding to quasi-Casimir $K_4=hkcfgS_4^2$. Note that, e.g., $K_1=edcbaS_1^2$. In the right picture we present the quiver after amalgamation.
	}
	\label{fi:Q2m}
\end{figure}
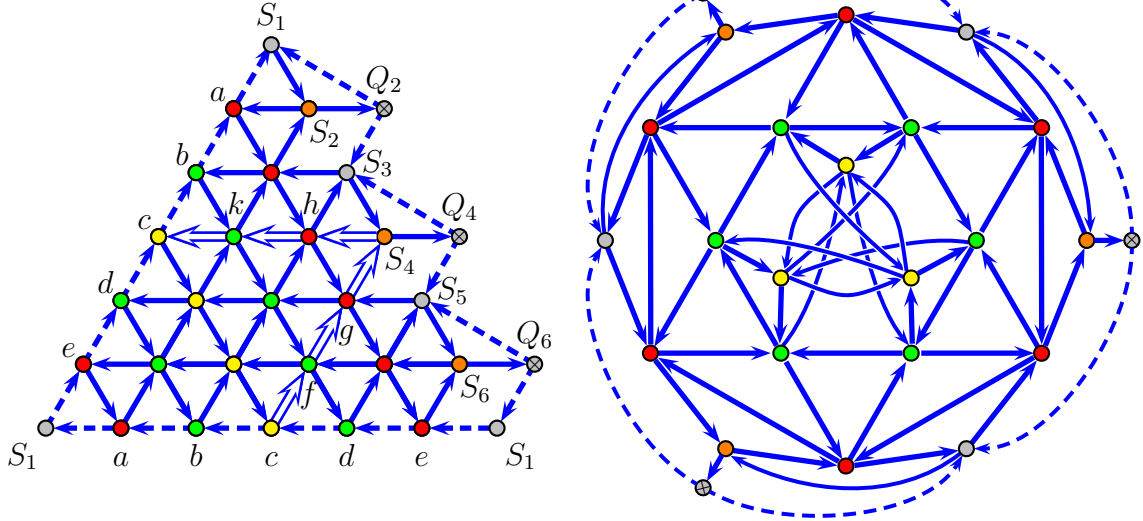

Constructing transport matrices $M_1$ and $M_2$ by the same rules as in Sec.~\ref{sec:QFG} and composing the block-upper triangular matrix $\mathbb A:=M_1^{\text{T}}M_2$ as in Sec.~\ref{s:A}, we have to amalgamate frozen variables on two sides of the $Sp_{2m}$ quiver thus obtaining the quiver depicted in Fig.~\ref{fi:Q2m}, which we call the  $\mathcal A_{Sp_{2m}}$-quiver.

Note that  $\mathcal A_{Sp_{2m}}$-quiver resembles  $\mathcal A_{n}$-quiver except for two rows on the right. We label variables of these rows by $S_l$ and $Q_{2r}$; each variable $S_l$ has a corresponding variable $Z_{(l,n-l,0)}$ in the $b_n$-quiver and $Q_{2r}$ are new variables.

We now find Casimirs for the  $\mathcal A_{Sp_{2m}}$-quiver. Note, first, that Casimirs $C_i$ (\ref{Ci}) remain Casimirs of  the  $\mathcal A_{Sp_{2m}}$-quiver. It was proved in \cite{ChM3} that for $mk\times mk$ block-triangular matrix $\mathbb A$ with blocks of size $k\times k$ we have $\bigl[\frac{mk}{2}\bigr]$ Casimirs $C_i$ and exactly $m\bigl[\frac{k+1}{2}\bigr]$ new Casimirs depending only on matrix entries $a_{i,j}$ inside diagonal $k\times k$ blocks. For $k=2$ we have exactly one such Casimir per each diagonal block, and this Casimir is just the determinant of the corresponding $2\times 2$ block. We now find the corresponding Casimirs for variables of $\mathcal A_{Sp_{2m}}$-quiver. 

We first introduce \emph{quasi-Casimirs} $K_l$, $l=1,\dots,2m$ having the same content as Casimirs $K_l$ (\ref{eq:Kl}) in $b_n$-quiver thus denoting by the same letter; we just replace $Z_{(l,n-l,0)}$ by $S_l$. An example of $K_4$ is in the left picture in Fig.~\ref{fi:Q2m}. Quasi-Casimirs $K_l$ commute with all cluster variables except $S_i$ and $Q_{2r}$. Substituting now $K_l$ for all $S_l$, we detach the part of the quiver containing only variables $K_l$ and $Q_{2r}$; this part looks as
$$
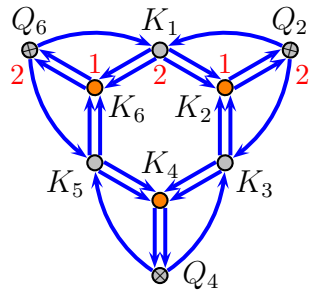
\begin{pspicture}(2,-2.3)(2,2.3)
\newcommand{\PATGEN}{%
{\psset{unit=1}
\psline[linecolor=blue,linewidth=1.5pt]{->}(-0.07,-1.1)(-0.07,-1.9)
\psline[linecolor=blue,linewidth=1.5pt]{->}(0.07,-1.1)(0.07,-1.9)
\psline[linecolor=blue,linewidth=1.5pt]{->}(-0.936,-0.4)(-0.936,0.4)
\psline[linecolor=blue,linewidth=1.5pt]{->}(-0.796,-0.4)(-0.796,0.4)
\psline[linecolor=blue,linewidth=1.5pt]{->}(0.936,-0.4)(0.936,0.4)
\psline[linecolor=blue,linewidth=1.5pt]{->}(0.796,-0.4)(0.796,0.4)
\psarc[linecolor=blue,linewidth=1.5pt]{<-}(0.866,-0.5){1.73}{184}{236}
\psarc[linecolor=blue,linewidth=1.5pt]{->}(-0.866,-0.5){1.73}{304}{356}
\pscircle[fillstyle=solid,fillcolor=lightgray](-0.866,-0.5){.1}
\pscircle[fillstyle=solid,fillcolor=orange](0,-1){.1}
\pscircle[fillstyle=solid,fillcolor=lightgray](0,-2){.1}
\rput(0,-2){\makebox(0,0)[cc]{\hbox{{\tiny$\times$}}}}
}}
\rput(0,0){\PATGEN}
\rput{120}(0,0){\PATGEN}
\rput{240}(0,0){\PATGEN}
\put(0,1.2){\makebox(0,0)[bc]{\hbox{{$K_1$}}}}
\put(0.7,0.4){\makebox(0,0)[tr]{\hbox{{$K_2$}}}}
\put(1,-0.6){\makebox(0,0)[tl]{\hbox{{$K_3$}}}}
\put(0,-0.7){\makebox(0,0)[bc]{\hbox{{$K_4$}}}}
\put(-1,-0.6){\makebox(0,0)[tr]{\hbox{{$K_5$}}}}
\put(-0.7,0.4){\makebox(0,0)[tl]{\hbox{{$K_6$}}}}
\put(-1.73,1.2){\makebox(0,0)[bc]{\hbox{{$Q_6$}}}}
\put(1.73,1.2){\makebox(0,0)[bc]{\hbox{{$Q_2$}}}}
\put(0.3,-2){\makebox(0,0)[cl]{\hbox{{$Q_4$}}}}
\put(0,0.8){\makebox(0,0)[tc]{\hbox{\tcr{\small $2$}}}}
\put(1.8,0.8){\makebox(0,0)[tl]{\hbox{\tcr{\small $2$}}}}
\put(-1.8,0.8){\makebox(0,0)[tr]{\hbox{\tcr{\small $2$}}}}
\put(-0.866,0.7){\makebox(0,0)[bc]{\hbox{\tcr{\small $1$}}}}
\put(0.866,0.7){\makebox(0,0)[bc]{\hbox{\tcr{\small $1$}}}}
\end{pspicture}
$$
It is easy to observe that, say, $Q_6^2K_6K_1^2K_2Q_2^2$ is a Casimir (numbers in the figure above indicate powers of the corresponding cluster variables). We therefore have the following statement.
\begin{lemma}\label{lem:Sp2m}
Casimirs of $\mathcal A_{Sp_{2m}}$-quiver are $m$ Casimirs $C_i$ (\ref{Ci}) and $m$ Casimirs $R_j$,
$$
R_j:=Q_{2j}^2 K_{2j} K_{2j+1}^2 K_{2j+2}Q_{2j+2}^2,
$$ 
with $K_l$ given by the same combinations (\ref{eq:Kl}) of cluster variables as in the case of $b_n$-quiver.
\end{lemma}

\section{Concluding remarks}\label{s:conclusion}
\setcounter{equation}{0}

In this paper, we have solved the problem of finding solutions to the eigenvalue problem $(\mathbb A-\lambda \mathbb A^\dagger)\psi=0$ for the quantum matrix $\mathbb A^\hbar$ from (\ref{A-quantum}). Our solution was based on the cluster-variable realization of entries of $\mathbb A^\hbar$ found in \cite{ChSh2}. We had also performed an extensive verification of our results  (in the semiclassical limit) using Mathematica computer code \cite{Wolfram}.

Note that Remark~\ref{rm:general} indicates that our results are by no means confined to a ``triangle'' $\Sigma_{0,1,3}$: they remain valid for any system undergoing the ``form dynamics'' described in the introduction. However, in the case of matrices $\mathbb A$ of a more general form, besides $[n/2]$ Casimirs generated by the same expression  $\det(\mathbb A-\lambda\mathbb A^{\text{T}})$ as in the upper-triangular case, we have additional Casimir operators related to ratios of minors located in the lower-left corner of the matrix $\mathbb A$ (see, e.g., \cite{ChM3}). We constructed the corresponding Casimirs for the  $\mathcal A_{Sp_{2m}}$-quiver corresponding to matrices $\mathbb A$ of block-upper triangular form with blocks of size $2\times 2$. This however does not change the results concerning the Jordan form of $\mathbb A\mathbb A^{-\dagger}$ and the corresponding symmetries because for every such system, $\mathbb A=M_1^{\text{T}}M_3 M_1$ with the same matrix $M_3$. Passing to more general $\mathbb A$ also does not change the conclusion that $\bigl[\mathbb A\mathbb A^{-\dagger}\bigr]^N=\mathbf I$ for any such matrix provided $\lambda_i$ constructed in Theorem~\ref{th:M} are distinct $N$th roots of unity.

\section*{Acknowledgements}
The authors are grateful to Alexander Shapiro for the useful discussion.
M.S. was supported by NSF grant DMS-1702115; both L.Ch. and M.S. are supported by International Laboratory of Cluster Geometry NRU HSE, RF Government grant, ag. \textnumero 2020-220-08-7077. H.Sh. thanks Department of Mathematics, Michigan State University and the exchange program "Discover America" for support and stimulating research atmosphere during his visit of MSU when the main part of the work was accomplished.


\begin{thebibliography}{99}

\footnotesize\itemsep=0pt





\bibitem{Bondal}
A.~Bondal, {\it A symplectic groupoid of triangular bilinear forms and the braid groups}, preprint IHES/M/00/02 (Jan. 2000).

\bibitem{BP99}
D. Bullock and J.H. Przytycki, {\it Multiplicative structure of Kauffman skein module quantization}, {\sl Proceedings Amer. Math. Soc.} {\bf 128(3)} (1999) 923--931










\bibitem{ChF2}
Chekhov L., Fock V., {\it Quantum mapping class group, pentagon relation, and geodesics},
{\sl Proc. Steklov Math. Inst.} {\bf 226} (1999), 149--163.

\bibitem{ChF3}
Chekhov L., Fock V., {\it Observables in 3d gravity and geodesic algebras},
{\sl Czech. J.~Phys.} {\bf 50} (2000) 1201--1208.


\bibitem{ChM}
Chekhov L.O., Mazzocco M., {\it Isomonodromic deformations and twisted Yangians arising in Teichm\"uller theory}, {\sl Advances Math.} {\bf 226(6)} (2011) 4731-4775, arXiv:0909.5350.



\bibitem{ChM3}
L. Chekhov and M. Mazzocco, {\it Block triangular bilinear forms and braid group action}, {\sl Comm. Math. Phys.} {\bf 322} (2013) 49--71.

\bibitem{ChM4}
L. Chekhov and M. Mazzocco, {\it On a Poisson space of bilinear forms with a Poisson Lie action}, {\sl Russ. Math. Surveys} {\bf 72(6)} (2017) 1109--1156; ArXiv:1404.0988v2.

\bibitem{CMR}
Chekhov L.O., Mazzocco M. and Rubtsov V., {\it Algebras of quantum monodromy data and decorated character varieties,}
{\it arXiv:1705.01447} (2017).






\bibitem{ChSh2}
L.O. Chekhov and M. Shapiro, {\it Darboux coordinates for symplectic groupoid and cluster algebras}, ArXiv:2003.07499v2, 43pp



\bibitem{Dub}
Dubrovin B., Geometry of $2$D topological field theories, Integrable systems and quantum groups (Montecatini Terme, 1993), {\it Lecture Notes in Math.,}\/ {\bf 1620}, Springer, Berlin, (1996) 120--348.









\bibitem{F97}
Fock V.V., {\it Dual Teichm\"uller spaces} arXiv:dg-ga/9702018v3, (1997).




\bibitem{FG1} V.~V.~Fock and A.~B.~Goncharov, {\it Moduli spaces of local systems and higher Teichm\"uller
theory}, {\sl Publ. Math. Inst. Hautes \'Etudes Sci.} 103 (2006), 1-211, math.AG/0311149 v4.





\bibitem{Fock-Rosly}
V.~V.~Fock and A.~A.~Rosly,
{\it Moduli space of flat connections as a Poisson manifold}, Advances in quantum field theory and statistical mechanics: 2nd Italian-Russian collaboration (Como, 1996), {\it Internat. J. Modern Phys. B} {\bf 11} (1997), no. 26-27, 3195--3206.




\bibitem{FST} 
Fomin S., Shapiro M., and Thurston D.,
{\it Cluster algebras and triangulated surfaces. Part I: Cluster complexes}, {\sl Acta Math.} {\bf 201} (2008), no. 1, 83--146.







\bibitem{GK91}
A.~M.~Gavrilik and A.~U.~Klimyk, {\it $q$-Deformed orthogonal and pseudo-orthogonal algebras
and their representations} {\sl Lett. Math. Phys.}, {\bf 21} (1991) 215--220.

\bibitem{GSV} M. Gekhtman, M. Shapiro, and A. Vainshtein, {\it Cluster algebra and Poisson geometry},
{\sl Moscow Math. J.} {\bf3}(3) (2003) 899-934.


\bibitem{Gold}
Goldman W.M., {\it Invariant functions on Lie groups and Hamiltonian
flows of surface group representations}, {\sl Invent. Math.} {\bf85}
(1986) 263--302.








\bibitem{HenTeit}
Henneaux~M. and Teitelboim~Cl., {\it Quantization of Gauge Systems}, 
{\sl Princeton University Press}, 1992.


\bibitem{Karasev}
M.~V.~Karasev, {\it Analogues of objects of Lie group theory by nonlinear Poisson brackets}, {\it Math. USSR Izvestia}
{\bf 28} (1987) 497--527.




\bibitem{KS}
Korotkin D. and Samtleben H.,
Quantization of coset space $\sigma$-models coupled to two-dimensional gravity,
{\it Comm. Math. Phys.} {\bf 190}(2) (1997) 411--457.


\bibitem{kirill}
Kirill Mackenzie 
General Theory of Lie Groupoids and Lie Algebroids,
{\it LMS Lect. Note Series} {\bf 213} (2005).




\bibitem{MRS} A. Molev, E. Ragoucy, P. Sorba, {\it Coideal subalgebras in quantum affine algebras}, {\sl Rev. Math. Phys.}, {\bf 15} (2003) 789--822.







\bibitem{NR} 
Nelson J.E., Regge T., Homotopy groups and $(2{+}1)$-dimensional
quantum gravity, {\it Nucl. Phys.~B} {\bf 328} (1989), 190--199.

\bibitem{NRZ}
Nelson J.E., Regge T., Zertuche F., Homotopy groups and
$(2+1)$-dimensional quantum de~Sitter gravity, {\it Nucl. Phys.~B} {\bf 339}
(1990), 516--532.






\bibitem{Po} Postnikov~A.,
Total positivity, Grassmannians, and networks, arXiv:math/0609764

\bibitem{PSW} Postnikov, A., Speyer, D., Williams, L., Matching polytopes, toric geometry, and the non-negative part of the Grassmannian, Journal of Algebraic Combinatorics 30 (2009), no. 2, 173-191. 


\bibitem{Gus-Al} G. Schrader and A. Shapiro, {\it Continuous tensor categories from quantum groups I: algebraic aspects}, arXiv:1708.08107.




\bibitem{Ugaglia}
Ugaglia M., {\it On a Poisson structure on the
space of Stokes matrices}, {\sl Int. Math. Res. Not.} {\bf 1999} (1999),  no.~9, 473--493,
{http://arxiv.org/abs/math.AG/9902045}{math.AG/9902045}.

\bibitem{Weinstein}
A. Weinstein, {\it Coisotropic calculus and Poisson groupoids}, {\sl J. Math. Soc. Japan} {\bf 40}(4) (1988) 705--727.

\bibitem{Wolfram}
S. Wolfram Mathematica, wolfram.com/mathematica/.

\bibitem{Ya} 
M. Yakimov, {\it Symplectic leaves of complex reductive Poisson--Lie groups}, {\sl Duke Math. J.} {\bf 112} (2002) 453--509.



\end{thebibliography}
\end{document}